%% file: AITgeo.tex
\documentclass[12pt]{article}

\usepackage{setspace}

\usepackage{amssymb}

\usepackage{amsmath}

\usepackage{mathrsfs}

\usepackage{stmaryrd}

\usepackage{enumitem}

\usepackage[all]{xy}

\usepackage{textcomp}

\usepackage{amsbsy}

\usepackage{epsfig,wrapfig}
\usepackage{graphicx}
\usepackage{array,multirow,booktabs}

\usepackage{url,comment}
\usepackage{color}

\input{Preamble_3}

\begin{document}

\title{Maximal linear spaces contained in the base loci of pencils of quadrics}
\author{Xiaoheng Wang}
\maketitle

\pagenumbering{roman}

\singlespacing

\input{abstract}

\tableofcontents\phantomsection

\pagenumbering{arabic}

\input{math}

\singlespacing

\addcontentsline{toc}{section}{Reference}
\bibliographystyle{abbrv}
\bibliography{biblio}

\end{document}

%% file: Preamble_3.tex
\usepackage{amscd,bbm}
\usepackage{amsfonts,mathrsfs,amssymb}
\usepackage{amsmath,amsthm,centernot}

\usepackage{url,comment}
\usepackage{color}
\usepackage[version=3]{mhchem}

\usepackage{epsfig,wrapfig}
\usepackage{graphicx}
\usepackage{caption}
\usepackage{subcaption}
\usepackage{array,multirow,booktabs}

\usepackage{varioref} 
\usepackage[all]{xy} 

\usepackage{setspace}
\usepackage{hyperref}
\hypersetup{colorlinks=true, linkcolor=black}

\theoremstyle{definition}
\newtheorem{theorem}{Theorem}[section]

\newtheorem{thm-defn}[theorem]{Theorem/Definition}
\newtheorem{lemma}[theorem]{Lemma}
\newtheorem{proposition}[theorem]{Proposition}
\newtheorem{cor}[theorem]{Corollary}

\theoremstyle{definition}

\newtheorem{example}[theorem]{Example}

\newtheorem{remark}[theorem]{Remark}

\numberwithin{equation}{section}

\def\be {{\mathbf e}}

\def\bba {{\mathbb A}}

\def\bbg {{\mathbb G}}

\def\bbp {{\mathbb P}}

\def\bbr {{\mathbb R}}

\def\bbz {{\mathbb Z}}

\def\scl {{\mathcal L}}

\def\sco {{\mathcal O}}

\def\ppU {\bbp(U)}
\def\pX {\bbp X}

\def\pY {\bbp Y}

\def\al {\alpha}
\def\be {\beta}

\def\vphi {\varphi}

\def\SO {\mbox{SO}}
\def\PSO {\mbox{PSO}}

\def\GL {\mbox{GL}}

\def\dim {\mbox{dim}\,}

\def\tDiv {\mbox{Div}}
\def\tdiv {\mbox{div}}

\def\tSpan {{\mbox{Span}}}

\def\tTr {{\mbox{Tr}}}

\def\tStab {\mbox{Stab}}

\def\PO {\mbox{PO}}
\newcommand{\Pic}[1] {\mbox{Pic}(#1)}

\def\Gal {\mbox{Gal}}
\def\picz {\mbox{Pic}^0}
\def\pico {\mbox{Pic}^1}
\def\dcup {\,\dot\cup}

\def\dim {\mbox{dim}}

\def\bpr {\vspace{-5pt}\pbpr \textbf{Proof: }}
\def\pbpr {}
\def\epr {\qed\vspace{10pt}}

\newcommand{\gl}[2]{{^{#2\!}#1}}

\def\to {\rightarrow}
\def\impl {\Rightarrow}

\def\onto {\twoheadrightarrow}
\def\disc {\mbox{disc}}
\def\Rlk {\mbox{Res}_{L/k}}

\renewcommand{\bar}[1]{\overline{#1}}
\newcommand{\w}[1]{\widetilde{#1}}
\newcommand{\lrg}[1] {\langle#1\rangle}

\newcommand{\udots}{\mathinner{\mskip1mu\raise1pt\vbox{\kern7pt\hbox{.}}\mskip2mu\raise4pt\hbox{.}\mskip2mu\raise7pt\hbox{.}\mskip1mu}}

\newcommand{\hide}[1]{{}}

\def\begEnu {\begin{enumerate}}
\def\begItem {\begin{itemize}}
\def\endEnu {\end{enumerate}}
\def\endItem {\end{itemize}}

%% file: abstract.tex

\centerline{Abstract}

\vspace{0.3in}

The geometry of the Fano scheme of maximal linear spaces contained in the base locus of a pencil of quadrics has been studied by algebraic geometers when the base field is algebraically closed. In this paper, we work over an arbitrary base field of characteristic not equal to 2 and show how these Fano schemes are related to the Jacobians of hyperelliptic curves. In particular, if $B$ is the base locus of a generic pencil of quadrics in $\bbp^{2n+1}$, and $F$ is the Fano variety of $n - 1$ planes contained in $B$, then $F$ is a component of a disconnected commutative algebraic group $G = \picz(C)\dcup F\dcup \pico(C)\dcup F'$, where $C$ is the hyperelliptic curve defined by the discriminant form of the pencil. In the second half of this paper, we study regular pencils of quadrics, where the hyperelliptic curve defined by the discriminant is singular.

%% file: math.tex

\section{Introduction}\label{sec:intro}

\input{Introduction}

\input{acknowledgements}


\section{Generic Pencil}\label{sec:generic}

\input{PenGeneric}


\section{Regular Pencil}\label{sec:regular}

\input{PenRegular}





%% file: Introduction.tex
Let $k$ be a field of characteristic not 2 and let $\scl = \{xQ_1-yQ_2|[x,y]\in \bbp^1\}$ be a pencil of quadrics in $\bbp^{N-1}$ for $N\geq3$ with $Q_1,Q_2$ defined over $k$. Let $B = Q_1\cap Q_2$ denote the base locus. In this paper, we study the general geometry of the variety $F$ of maximal dimensional linear subspaces contained in $B$. Let $A_1,A_2\in M_N(k)$ denote two Gram matrices of $Q_1,Q_2$ respectively and let $f(x)$ be the polynomial of degree at most $N$ defined by
$$f(x) = \disc(xA_1 - A_2) = (-1)^{N(N-1)/2}\det(xA_1 - A_2).$$
We shall assume that $f(x)$ splits completely over a separable closure $k^s$ of $k$..

The geometry depends very much on the parity of $N$ and how ``singular'' the pencil is. A pencil is \textbf{generic} if $\scl$ is a generic line in the space $\bbp(H^0(\sco_{\bbp^{N-1}}(2)))$ of all quadrics. Equivalently, $\scl$ contains precisely $N$ singular quadrics over $k^s$ which are all simple cones. This is also equivalent to saying that $f(x)$ has degree at least $N-1$ with no repeated roots. Suppose $\scl$ is generic, denote by $C$ the hyperelliptic curve with affine equation $y^2 = f(x)$. The isomorphism type of $C$ over $k$ is independent of the choices of $A_1,A_2$.

When $N=2n+1$ is odd, the dimension of maximal linear subspaces contained in $B$ is $n-1$. We assume that $k$ has at least $N+1$ elements which implies that there is a rational non-degenerate quadric over $k$ in the pencil. After renaming, we may assume that $Q_1$ is non-degenerate and hence $C$ has genus $n$. It is well-known from intersection theory that geometrically over $k^s$, there are $2^{2n}$ such $n-1$ planes. The arithmetic theory over $k$ was studied in \cite{BG} where it was shown that $F$ is a torsor of $J[2]$, where $J$ denotes the Jacobian of $C$. When $N=2n+2$ is even, the theory is richer. The dimension of maximal linear subspaces contained in $B$ is still $n-1$ and $C$ has genus $n$. The rational function $x$ defines a degree 2 morphism $C\rightarrow \bbp^1.$ Let $D_0$ denote the hyperelliptic class obtained from pulling back the hyperplane section on $\bbp^1.$ It was proved by Desale and Ramanan \cite{DR}, Reid \cite{Reid}, and Donagi \cite{Donagi} that geometrically over $k^s$, $F$ is isomorphic to the Jacobian $J$ of $C$. As Weil pointed out in \cite{Weil}, Gauthier had studied this in \cite{Gauthier}. The arithmetic theory when $C$ has genus 1 is known and is used heavily in studying the 4-descent of elliptic curves. The main result of this paper is that for arbitrary $n\geq1,$ the Fano scheme $F$ is a torsor of $J$ over $k$ and moreover,

\begin{theorem}\label{thm:Main}
There is a commutative algebraic group structure $+_G$ over $k$ on the disconnected variety $$G = \picz(C) \dcup F \dcup \pico(C) \dcup F'$$ such that,
\begEnu
\item $G^0 = \picz(C)$ with component group $G/G^0 \simeq \bbz/4,$
\item $F'$ is isomorphic to $F$ as varieties via the inversion map $-1_G$,
\item the group law extends that on $H = \Pic{C}/D_0 \simeq \picz(C) \dcup \pico(C)$.
\endEnu
\end{theorem}

Moreover, we will show that this structure is unique once we impose one more condition. See Theorem \ref{thm:maintheorem} for the complete statement. Note having such a group structure is strictly better than knowing only that $2[F]=[\pico(C)]$ for it also gives a canonical lift of $[F]$ to a torsor of $J[4]$ by taking
$$F[4] = \{X\in F|X +_G X +_G X +_G X = 0\}.$$
When $\pico(C)(k)\neq\emptyset$, $F$ is a torsor of $J$ of order dividing 2 and for any $[D_1]\in\pico(C)(k),$ one has a canonical lift of $[F]$ to a torsor of $J[2]$ by taking
$$F[2]_{[D_1]} = \{X\in F| X +_G X = [D_1]\}.$$
When the class $[D_1]$ comes from a rational point $P$ on $C$, the lift $F[2]_P$ can also be described geometrically, see Example \ref{ex:Weier} and \ref{ex:nonWeier}. These two special cases are the key geometric input used in \cite{BG3} and \cite{SW6} to obtain the average sizes of the 2-Selmer groups of Jacobians of hyperelliptic curve with a rational Weierstrass point, or a rational non-Weierstrass point, respectively.

In the second half of the paper, we generalize these results to the case where the pencil $\scl$ is ``slightly singular'', or \textbf{regular}. A pencil $\scl$ is regular if it only has simple cones as singular members. In other words, $\scl$ could intersect tangentially to the discriminant hypersurface in $\bbp(H^0(\sco_{\bbp^{N-1}}(2)))$ but does not contain any quadrics with higher degeneracy degree than simple cones. Let $U$ denote the underlying $k$-vector space of dimension $N$ and view the two generating quadrics $Q_1,Q_2$ as linear operators $A_1,A_2:U\rightarrow U^*.$ When $Q_1$ is smooth, $A_1$ is an isomorphism and the composite
$$T:U\xrightarrow{A_2}U^*\xrightarrow{A_1^{-1}}U$$
is self-adjoint with respect to $A_1$. We call $T$ the self-adjoint operator associated to the pencil. The pencil spanned by $Q_1,Q_2$ is regular if and only if $T$ is regular, which by definition means that all the eigenspaces of $T$ are 1-dimensional. Note the pencil is generic if and only if $T$ is regular semi-simple.

Suppose $\scl$ is regular with $Q_1$ non-degenerate and $N=2n+1$ is odd. One may then assume $f(x)$ is monic and it factors as $f(x) = \prod_{i=1}^{r+1}(x - \al_i)^{m_i}$ over $k^s$. Let $U_{i,T}$ denote the generalized eigenspace over $k^s$ of $T$ with eigenvalue $\al_i.$ Since $T$ is self-adjoint with respect to $A_1$, its generalized eigenspaces are pairwise orthogonal. A projective $n-1$ plane contained in the base locus $B$ can be viewed as a linear $n$ plane $X$ such that $X\subset X^\perp,TX\subset X^\perp$ where $\perp$ is taken with respect to $A_1$. For each $i=1,\ldots,r+1,$ we define $\dim_{i,T}(X)$ to be the dimension of the maximal $T$-stable subspace of $(X\otimes k^s)\cap U_{i,T}.$ Since each $U_{i,T}$ is $m_i$ dimensional and $A_1$ restricts to a non-degenerate quadratic form on $U_{i,T}$, we have $$\dim_{i,T}(X) \leq m_i/2.$$
For any sequence of integers $d_1,\ldots,d_{r+1}$ such that $0\leq d_i\leq m_i/2,$ we define
$$L^T_{\{d_1,\ldots,d_{r+1}\}}(k^s)=\{X\simeq (k^s)^n|X\subset X^\perp,TX\subset X^\perp, \dim_{i,T}(X)= d_i\}.$$
Note the singular locus of $B$ consists of the projectivization of the eigenspaces of $T$ whose associated eigenvalue has multiplicity at least 2. Hence $L^T_{\{0,0,\ldots,0\}}$ is the set of $k^s$-points of the variety of projective $n-1$ planes contained in the smooth part of $B$.

\begin{theorem}\label{thm:MainRegOdd}
$|L^T_{\{d_1,\ldots,d_{r+1}\}}(k^s)| = 2^r/2^a$ where $a$ is the number of $d_i$'s equal to $m_i/2.$
\end{theorem}

Now for any field $k'$ containing $k$, one defines $L^T_{\{d_1,\ldots,d_{r+1}\}}(k')$ similarly. Let $J=\picz(C)$ denote the generalized Jacobian of the (possibly singular) complete curve $C$ defined by affine equation $y^2 = f(x)$. Then,

\begin{theorem}\label{thm:MainRegOddAIT}
For any field $k'$ containing $k$, $J[2](k')$ acts on $L^T_{\{d_1,\ldots,d_{r+1}\}}(k')$ simply-transitively if $a=0$ and transitively if $a>0$.
\end{theorem}

As before, the case when $N=2n+2$ is even is more interesting. In what follows, we will use $X$ as a linear subspace of $U$, $\pX$ its projectivization. For any $v\in U\otimes k^s$, denote by $[v]$ the point of $\bbp(U)(k^s)$ corresponding to the line spanned by $v$. As before, let $C$ denote the hyperelliptic curve defined by affine equation $y^2 = f(x).$ Then $C$ is smooth if and only if the pencil is generic. When all the roots of $f$ have multiplicity at most 2, $C$ is nodal. Denote by $p_g$ the geometric genus of $C$, defined as the genus of the normalization $\w{C}$. When $p_g=-1,$ the curve $C$ is reducible and the base locus $B$ contains one $\bbp^{n}.$ When $p_g\geq0,$ $B$ contains no $\bbp^n$ and we define $$F_0 = \{\pX|\dim\pX = n-1,\pX \subset B\}.$$
As in the odd case, we need to impose certain open conditions to obtain interesting relations with the Jacobian of $C$. 
Consider
$$F=\{\pX\in F_0|[v]\notin\pX,\mbox{ for all singular points }[v]\in B\}.$$
Bhosle \cite{Bhosle} proved that when $C$ only has nodal singularities, $F$ is isomorphic to the generalized Jacobian of $C$ over $k^s$. Over the base field $k$, we have the following result.

\begin{theorem}\label{thm:MainRegEven}
Suppose $p_g\geq0$ and $C$ only has nodal singularities. Then there is a commutative algebraic group structure $+_G$ defined over $k$ on the disconnected variety $$G = \picz(C) \dcup F \dcup \pico(C) \dcup F'$$ such that,
\begEnu
\item $G^0 = \picz(C)$ with component group $G/G^0 \simeq \bbz/4,$
\item $F'$ is isomorphic to $F$ as varieties via the inversion map $-1_G$,
\item the group law extends that on $H = \Pic{C}/D_0 \simeq \picz(C) \dcup \pico(C)$ where $D_0$ is the hyperelliptic class.
\endEnu
\end{theorem}

\begin{cor}\label{cor:compactification}
Over $k^s$, one can take $F_0$ as a compactification of the generalized Jacobian $\picz(C).$ 
\end{cor}

We expect the above theorem to be true without the condition on $C$. Let $\w{F}$ denote the torsor of $J(\w{C})$ obtained from certain reductions on the pencil $\scl$ to the generic case. We can prove, without the condition on $C$, that,

\begin{theorem}\label{thm:MainRegEvenCarb}
There is a surjection $F\rightarrow\w{F}$. Over $k^s$, the pre-image of every point has a filtration with $\bbg_a$ and $\bbg_m$ factors. The kernel of the natural map $J(C)\to J(\w{C})$ has a filtration with the same factors.
\end{theorem}

%% file: acknowledgements.tex

\section*{Acknowledgements}

This paper is a revised version of part of my Harvard Ph.D. thesis written under the supervision of Benedict H. Gross. I am very grateful to him for introducing me to the subject and for many insightful discussions. I would also like to thank Manjul Bhargava, George Boxer, Brian Conrad, Anand Deopurkar, Joe Harris and Anand Patel for valuable suggestions and conversations.





%% file: PenGeneric.tex
Let $k$ be a field of characteristic not 2 and let $Q_1,Q_2$ be two quadratic forms on a $k$-vector space $U$ of dimension $N$. In this chapter, we study the general geometry of the maximal isotropic subspaces with respect to both quadrics.

There are three equivalent ways to formulate this problem. We call the above formulation the $(Q_1,Q_2)-$setup. Suppose now $Q_1$ is non-degenerate. Let $b_1,b_2$ denote the corresponding bilinear form, $$b_i(v,w) = \frac{1}{2}(Q_i(v+w) - Q_i(v) - Q_i(w)).$$ Let $T:U\rightarrow U$ be the unique operator such that for all $v,w\in U,$
\begin{equation}\label{eq:defT}
b_2(v,w)=b_1(v,Tw).
\end{equation}
Note $T$ is self-adjoint with respect to $b_1$ since $b_1,b_2$ are symmetric.

To say a linear subspace $X$ is isotropic with respect to both $Q_1,Q_2$ is the same as saying
\begin{equation}\label{eq:XTXp}
X\subset X^{\perp_{Q_1}}, TX\subset X^{\perp_{Q_1}}.
\end{equation}
Therefore, instead of starting with a pair of quadratic forms, we could have started with a non-degenerate quadratic form along with a self-adjoint operator. We call this formulation the $(Q_1,T)-$setup.

Lastly, we could view $Q_1,Q_2$ as quadrics in $\ppU$ and take a pencil $\scl = \{xQ_1-yQ_2|[x,y]\in \bbp^1\}$ of quadrics in $\ppU.$ Let $B=Q_1\cap Q_2$ denote the base locus. The above problem regarding common isotropic subspaces translates into studying the Fano variety of linear subspaces contained in the base locus. We call this formulation the $(\ppU,\scl)-$setup.

We define the notion of \textbf{generic} in each of the three formulations. With the $(Q_1,Q_2)$-setup, we require $f(x)=(-1)^{N(N-1)/2}\det(xA_1-A_2)$ to have no repeated roots where $A_1,A_2$ are two Gram matrices for $Q_1,Q_2$ respectively. With the $(Q,T)$-setup, we require the characteristic polynomial $f_T(x)=\det(xI-T)$ of the self-adjoint operator $T$ to have no repeated roots. We will also assume that $k$ has at least $N+1$ elements, for otherwise there might not exists a rational non-degenerate $Q$ in the pencil. With the $(\ppU,\scl)$-setup, we require that the pencil $\scl$ is a generic line in $\bbp(H^0(\sco_{\bbp^{N-1}}(2)))$. Equivalently, $\scl$ contains precisely $N$ singular quadrics over $k^s$ which are all simple cones.

\subsection{Odd dimension}\label{sec:povnonsin}

Suppose $U$ has dimension $N=2n+1$ and $Q_1$ is non-degenerate. Let $C$ be the hyperelliptic curve of genus $n$ defined by affine equation
$$y^2 = f(x) = (-1)^n\det(xA_1 - A_2).$$
The isomorphism type of $C$ over $k$ is independent of the choices of the Gram matrices $A_1,A_2$. Let $J$ denote the Jacobian of $C$. We assume that $f(x)$ splits completely over $k^s$ and adopt the $(Q,T)$-setup for this case. For every field $k'$ containing $k$, let $W_T(k')$ denote the set of (linear) $n$-dimensional $k'$-subspaces $X$ of $U\otimes k'$ such that $X\subset X^\perp, TX\subset X^\perp.$

The geometry over $k^s$ is well-known using classical intersection theory. The following explicit description of the $2^{2n}$ elements of $W_T(k^s)$ is due to Elkies. After a change of basis over $k^s$, one can assume the two quadrics are given by,
\begin{eqnarray*}
Q_1(x) &=& x_1^2 + x_2^2 + \cdots + x_{2n+1}^2\\
Q_2(x) &=& c_1x_1^2 + c_2x_2^2 + \cdots + c_{2n+1}x_{2n+1}^2.
\end{eqnarray*}
The following system of linear equations in $D_1,\ldots,D_{2n+1}$ has a 1-dimensional kernel.
\begin{eqnarray*}
D_1 + D_2 + \cdots + D_{2n+1} &=& 0\\
c_1D_1 + c_2D_2 + \cdots + c_{2n+1}D_{2n+1} &=& 0\\
&\vdots&\\
c_1^{2n-1}D_1 + c_2^{2n-1}D_2 + \cdots + c_{2n+1}^{2n-1}D_{2n+1} &=& 0
\end{eqnarray*}
Choose a basis $(D_1,\ldots,D_{2n+1})$ for the kernel and note that $D_i\neq 0$ for all $i$. For each $i=1,\ldots,2n+1$, let $d_i$ be a square root of $D_i.$ Any choice of a system of square roots gives an element of $W_T(k^s)$ by taking
$$X = \{(d_1P(c_1),\ldots,d_{2n+1}P(c_{2n+1}))|P\mbox{ any polynomial of degree }n-1\}.$$

The arithmetic aspect of the theory has been studied by Bhargava and Gross in \cite{BG}, we list the results without proof. Consider the following two schemes over $k$.
$$V_f = \{T:U\rightarrow U | T^* = T\mbox{ with characteristic polynomial }f\}\subset \bba^{(2n+1)^2},$$
$$W_f = \{(T,X)\in V_f\times Gr(n,U)|X\subset X^\perp, TX\subset X^\perp\}.$$
Note $W_T$ is the fiber of $W_f$ above a fixed $T$. The group $\PO(U,Q)=O(U,Q)/(\pm1)$ acts on $V_f, W_f$ via
$$g.T = gTg^{-1},\quad g.(T,X) = (gTg^{-1},gX).$$

\begin{proposition}\label{prop:POtransitive}
The action of $\PO(U,Q)$ on $V_f$ has a unique geometric orbit. For any $T\in V_f(k')$ defined over some field $k'$ over $k$, its stabilizer scheme $\tStab(T)$ is isomorphic to $\mbox{Res}_{L'/k'}\mu_2/\mu_2\simeq J_{k'}[2]$ as group schemes over $k'$ where $L' = k'[x]/f(x).$
\end{proposition}

For general $Q$, there might not be a self-adjoint operator defined over $k$ with the prescribed characteristic polynomial. For example over $\bbr$, operators self-adjoint with respect to the positive definite form have real eigenvalues.

\begin{lemma}\label{lem:split} If $Q$ is split, then $V_f(k)$ and $W_f(k)$ are nonempty. Furthermore, there exists $(T_0,X_0)\in W_f(k)$ with trivial stabilizer in $\PO(U,Q)(k^a).$
\end{lemma}

\begin{theorem}\label{thm:simpletransitivesmoothclose} Suppose $k$ is algebraically closed or separably closed, then $\PO(V,Q)(k)$ acts simply-transitively on $W_f(k).$
\end{theorem}

\begin{cor}\label{thm:simpletransitivesmooth} Suppose $k$ is arbitrary. Then $\PO(V,Q)(k')$ acts simply-transitively on $W_f(k')$ for any field $k'$ over $k$.
\end{cor}

\begin{cor}\label{cor:StabJ2} Suppose $k$ is arbitrary, and $T\in V_f(k).$ Let $J$ denote the Jacobian of the hyperelliptic curve defined by $y^2 = f(x).$ Then $W_T$ is a torsor for $J[2].$
\end{cor}

\begin{remark}\label{rmk:J2inPO}
One can write down an explicit formula for the identification
\begin{equation}\label{eq:StabJ2}
J[2]\simeq\tStab(T).
\end{equation}
We will work over $k^s$ and it will be clear that the map is Galois equivariant. Denote the roots of $f(x)$ over $k^s$ by $\al_1,\ldots,\al_{2n+1},$ and by $P_i$ the Weierstrass point corresponding to the root $\al_i.$ Recall $J[2]$ is an elementary 2-group generated by $(P_i)-(\infty)$ with the only relation being that their sum is trivial. For each generator $(P_i)-(\infty)$, one looks for a polynomial $g_i(x)$ such that $g_i(\al_i)=-1$ and $g_i(\al_j)=1$ for all $j\neq i.$ Then $g_i(T)$ is the image of $(P_i)-(\infty)$ in $\tStab(T).$ The image does not depend on the choice of the polynomial $g_i$ because any two choices defer by some multiples of $f(x)$ and $f(T)=0.$ Define $h_i(x) = f(x)/(x-\al_i),$ then $$g_i(x) = 1 - 2\frac{h_i(T)}{h_i(\al_i)}$$ does the job. In other words, on the level of $k^s$-points, \eqref{eq:StabJ2} is given by
$$\sum (\al_i)-(\infty) \mapsto \prod\left(1-2\frac{h_i(T)}{h_i(\al_i)}\right) = 1 - 2\sum \frac{h_i(T)}{h_i(\al_i)}.$$
The above summation and product are written without indices, meaning the above equality holds for any (finite) collection of matching indices.

It is also worth noting that the commutativity of $\tStab(T)$ is necessary for this identification to be canonical. Suppose $T'$ is in the same $\PO$-orbit as $T$, and let $g$ be an element in $G$ such that $gTg^{-1}=T'.$ The composite map $$\tStab(T)\rightarrow J[2]\rightarrow \tStab(T')$$ is given by conjugation by $h\mapsto ghg^{-1}$. Naturality requires this map to be independent on the choice of $g$. Now if $g'$ also sends $T$ to $T'$, then there exists some $g_0\in \tStab(T)$ such that $g'=gg_0.$ Conjugation by $g'$ and $g$ induce the same map on $\tStab(T)$ because conjugation by $g_0$ is trivial due to the commutativity of $\tStab(T).$

See Remark \ref{rmk:J2inPOalt} for a different view point of \eqref{eq:StabJ2}.
\end{remark}

\subsection{Even dimension}\label{sec:torofJ}

Suppose $U$ has dimension $2n+2$. The projective formulation is easier to work with in this case.

Let $\scl = \{Q_\lambda|\lambda\in \bbp^1\}$ be a rational generic pencil of quadrics in $\bbp^{2n+1}=\ppU.$ Rationality means it is generated by two quadrics $Q_1,Q_2$ defined over $k$.

The cone points of the $2n+2$ singular quadrics are best understood in terms of the self-adjoint operator $T$ defined in \eqref{eq:defT} assuming $Q_1$ is non-degenerate. The quadric $\lambda Q_1 - Q_2$ is singular if and only if $\lambda$ is an eigenvalue of $T$. If we denote a corresponding eigenvector by $v_\lambda,$ then the cone point of $\lambda Q_1 - Q_2$ is $[v_\lambda]\in\ppU.$ In particular, the $2n+2$ cone points span the entire $\ppU$.

Since $\scl$ is generic, the maximal (projective) dimension of any linear space contained in the base locus $B$ is $n-1$. Consider the following variety over $k$,
$$F = \{\pX | \dim(\pX)=n-1, \pX\subset B\}.$$

\subsubsection*{The hyperelliptic curve $C$}

For any rational generic pencil $\scl,$ there is an \emph{associated hyperelliptic curve} defined as follows.

For any quadric $Q$ in $\bbp^{2n+1}$, one defines its Lagrangian variety by $$L_Q = \{\pY | \pY\subset Q, \dim(\pY) = n\}\subset\bbg r(n, \ppU).$$ When $Q$ is smooth, $L_Q$ has two connected components, also called the \textbf{rulings} of $n$-planes in $Q$. Two $n$-planes in $Q$ lie in the same ruling if and only if their intersection codimension in either one of them is even. If $Q$ is defined over some field $k',$ its \textbf{discriminant} is defined by $$\disc(Q) = (-1)^{n+1}\det(Q)k'^{\times2}\in k'^\times/k'^{\times2}.$$ The connected components of $L_Q$ are defined over $k'(\sqrt{\disc(Q)}).$ In other words, $L_Q(k'^s)$ hits both rulings and the $\Gal(k'^s/k'(\sqrt{\disc(Q)}))$-action on $L_Q(k'^s)$ preserves the rulings. When $Q$ is singular, $L_Q$ has only one connected component.

Consider the following variety
$$\w{F} = \{(Q_\lambda, \pY)|\lambda\in\bbp^1, \pY\in L_{Q_\lambda}\}\subset \scl\times\bbg r(n, \ppU).$$
There is an obvious projection map $p_1:\w{F}\rightarrow\bbp^1.$ The fiber over $\lambda\in \bbp^1$ is isomorphic to $L_{Q_\lambda}.$ Let $$\epsilon:\w{F}\rightarrow C,\quad \pi:C\rightarrow\bbp^1$$ denote the Stein factorization. In other words, $\epsilon$ has connected fibers and the fibers of $\pi$ correspond bijectively to the connected components of the fibers of $p_1.$ Therefore, $C$ is a double cover of $\bbp^1$ branched over the $2n+2$ points that correspond to the singular quadrics on the pencil. A homogeneity analysis as in \cite{Donagi} shows that $C$ is smooth at the ramification points. Hence $C$ is a hyperelliptic curve of genus $n$, and to give a point on $C$ is the same as giving a quadric on the pencil plus a choice of ruling. We call $C$ the hyperelliptic curve associated to the pencil and it parameterizes the rulings in the pencil. The Weierstrass points of $C$ correspond to the $2n+2$ points on $\bbp^1$ cut out by $\det(xA_1 - yA_2).$ The curve $C$ is isomorphic over $k$, but not canonically, to the hyperelliptic curve defined by the affine equation
$$y^2 = (-1)^{n+1}\det(xA_1 - A_2).$$

It was known to algebraic geometers (\cite{Reid}, \cite{DR}, \cite{Donagi}) that when $k$ is algebraically closed, $F$ is isomorphic to $J$, the Jacobian of the curve $C$ defined above. Therefore it is natural to expect that over a general field, $F$ is a torsor of $J$. In fact, we prove something stronger:

\begin{theorem}\label{thm:Main2}
There is a commutative algebraic group structure $+_G$ over $k$ on the disconnected variety $$G = \picz(C) \dcup F \dcup \pico(C) \dcup F'$$ such that,
\begEnu
\item $G^0 = \picz(C)$ with component group $G/G^0 \simeq \bbz/4,$
\item $F'$ is isomorphic to $F$ as varieties via the inversion map $-1_G$,
\item the group law extends that on $H = \Pic{C}/D_0 \simeq \picz(C) \dcup \pico(C)$ where $D_0$ is the hyperelliptic class.
\endEnu
\end{theorem}

Moreover, we will show that this structure is unique once we impose one more condition. See Theorem \ref{thm:maintheorem} for the complete statement.

\subsubsection*{The dimension of $F$}

Since $F$ is isomorphic to $J$ over $k^s$, one can conclude that $F$ has dimension $n$ as an algebraic variety. Even without passing to the separable closure, one can still show that $F$ has dimension at least $n$. For any quadric $Q$ in $\bbp^{2n+1}$, let $F_{n-1,Q}$ denote the variety of $(n-1)$-planes in $Q$. When $Q$ is smooth, $F_{n-1,Q}$ is smooth irreducible of dimension $n(n+3)/2$ (\cite{Harris} p.293). Let $Q,Q'$ be two smooth quadrics on the pencil, then $F = F_{n-1,Q}\cap F_{n-1,Q'}$ has dimension at least
$$n(n+3)/2 + n(n+3)/3 - \mbox{dim}\bbg r(n-1,2n+1) = n.$$


\subsubsection*{A morphism $\tau:C\times F\rightarrow F$ that will serve as subtraction}

Given any pair $(c,\pX)\in C\times F$, there is a unique $n$-plane $\pY$ containing $\pX$ in the ruling of the quadric defined by $c$. This defines a morphism $$\delta:C\times F\rightarrow \w{F}.$$ The graph of $\delta$ is the following closed subvariety of $C\times F\times \w{F},$
$$\Sigma = \{(c,\pX,Q_{\pi(c)},\pY)|\pY\in L_{Q_{\pi(c)}}\mbox{ is in the ruling of }c,\pX\subset\pY\}.$$

Given any $(Q_{\lambda},\pY)\in\w{F}$ and $\pX\in F$ such that $\pX\subset\pY,$ let $Q$ be another quadric on the pencil. Since the base locus contains no $n$-planes, $\pY\cap B = \pY\cap Q$ is a quadric in $\pY$ containing $\pX$. Hence, $\pY\cap B=\pX\cup \pX'$ is the union of two (possibly equal) $(n-1)$-planes. We define
\begin{equation}
\tau(c,\pX)=\pX'.
\end{equation}

For a fixed $c\in C,$ define $\tau(c):F\rightarrow F$ by $\tau(c)\pX = \tau(c,\pX).$ Note $\tau(c)$ is an involution in the sense $\tau(c)^2 = \text{id}.$

One can write down a more explicit formula for $r$ as follows. Given any $(\pY,\pX)\in \Sigma,$ since $\pY\nsubseteq Q,$ there exists $p\in Y\backslash X$ such that $b(p,p)\neq0$ where $b$ is the bilinear form associated to $Q$. There is a linear map on $U\otimes k^a$ given by reflection about $p^{\perp_Q},$ namely
\begin{equation}\label{eq:defoftau}
\mbox{refl}_p:v\mapsto v - 2\frac{b(v,p)}{b(p,p)}p.
\end{equation}
Then $$r(\pY,\pX) = \bbp(\mbox{refl}_p(X)).$$

In order to put a group structure on $G = J\dcup F\dcup \pico(C)\dcup F',$ it suffices to define a simply-transitive action of $H = \Pic{C}/D_0$ on $F\dcup F'$ for then one can define $+_G$ as follows: for any $x,x'\in F\dcup F',[D],[D']\in H:$
\begEnu
\item $[D]+_G[D']$ is the usual addition in $H$,
\item $x +_G [D]=x + [D]$ is the image of $x$ under the action of $[D]$,
\item $x +_G x'$ is the unique element in $H$ that sends $-x'$ to $x$.
\endEnu

\subsubsection*{An action of $\tDiv(C)$ on $F\dcup F'$}

We start from the following action of $C$ on $F\dcup F':$
\begin{equation}\label{eq:actionofC}
\pX + (c) = -\tau(\bar{c})\pX,\quad\quad -\pX + (c) = \tau(c)\pX,
\end{equation}
where $c\mapsto\bar{c}$ denotes the hyperelliptic involution. The second equality follows the idea that $\tau:C\times F \rightarrow F$ serves as a subtraction, and the first equality was rigged so that divisors linearly equivalent to the hyperellipitic class $D_0$ acts trivially. The following Lemma allows one to extend this action to the semi-group of effective divisors on $C$. Negating \eqref{eq:actionofC} then gives the extension to the entire group of divisors.

\begin{lemma}\label{lem:Xe1e2} For any $x\in F\dcup F',c_1,c_2\in C,$
$$(x + (c_1)) + (c_2) = (x + (c_2)) + (c_1).$$
\end{lemma}

\bpr Unwinding the above definition, we need to prove for any $\pX\in F,$
$$\tau(c_2)\tau(\bar{c_1})\pX = \tau(c_1)\tau(\bar{c_2})\pX.$$
As both sides are defined by polynomial equations, it suffices to prove this equality for generic $\pX,c_1,c_2,$ over the algebraic closure, in particular we may assume there is no tangency involved. This is proved in \cite{Donagi} p.232 by looking at the following intersection $$\tSpan\{\pX,\tau(\bar{c_1})\pX,\tau(\bar{c_2})\pX\}\bigcap B.\epr$$

\begin{theorem}\label{thm:Hsimpletransitive} The above action of $\tDiv(C)$ descents to a simply-transitive action of $H$ on $F\dcup F'.$
\end{theorem}

\begin{remark}
\begEnu
\item $+_G$ is defined over $k$, because $\tau$ is defined over $k$.
\item $+_G$ is commutative. If $[D]$ sends $-x'$ to $x,$ it also sends $-x$ to $x'.$ This follows from the definition of the action of $\tDiv(C)$ on $F\dcup F'.$
\endEnu
\end{remark}

Before proving this Theorem, we give some concrete examples of $+_G$ in certain simple cases.

\begin{example}\label{ex:n=1}
Suppose $n=1$, then $F$ is the variety of points in the intersection of two generic quadrics in $\bbp^3$ and $C$ is a genus 1 curve. Given two points $\pX,\pX'\in F,$ let $\pY$ denote the line passing through them. There exists a unique quadric in the pencil and a unique ruling that contains $\pY,$ and this data is equivalent to giving a point on $C$. If one passes to the algebraic closure and identify $F\simeq J\simeq C,$ then $+_G:F\times F\rightarrow \pico(C)$ is just the addition on $J$.
\end{example}

\begin{example}\label{ex:codim1}
Suppose now $n$ is general and $\pX,\pX'\in F$ intersect in codimension 1 in either/both of them. Let $\pY = \tSpan\{\pX,\pX'\}$ denote their linear span, then $\pY\simeq\bbp^n.$ Let $p$ be a point on $\pY\backslash(\pX\cup\pX').$ There is a quadric $Q$ in $\scl$ containing $p$. Its intersection with $\pY$ contains two $\bbp^{n-1}$ and a point not on them, hence it cannot be a quadric. Furthermore, since the pencil is generic, the base locus contains no $\bbp^n$. Therefore, $\pY$ is contained in a unique quadric $Q$ in $\scl$ and a unique ruling on $Q$. Once again, such data determines a point on $c\in C$ and our group law says $$\pX +_G \pX' = (c)\in\pico(C).$$
\end{example}

\begin{example}\label{ex:bigTxB}
For any $\pX\in F$, since $B$ is a complete intersection, $$T_{\pX}B=T_{\pX}Q_1 \cap T_{\pX}Q_2 = \bbp(X^{\perp_{Q_1}} \cap X^{\perp_{Q_2}}).$$ As the next Lemma shows, $T_{\pX}B$ has dimension at most $n$. If $\pX\in F$ such that $T_{\pX}B\simeq \bbp^n,$ then just as in the above example, there exist a unique quadric in $\scl$ and a unique ruling that contains $T_{\pX}B.$ Such data determines a point on $c\in C$ and our group law says $$\pX +_G \pX = (c)\in\pico(C).$$ As we will see in Example \ref{ex:Weier}, for each Weierstrass point, there exists $2^{2n}$ such $\pX$ for which $T_{\pX}B\simeq\bbp^n$ is contained in the corresponding singular quadric.

\begin{lemma}\label{lem:dimTxB}
For generic pencil $\scl,$ $\dim(T_{\pX} B)\leq n.$
\end{lemma}

\bpr Suppose without loss of generality $Q_1,Q_2$ are non-degenerate. Since $\dim(X) = n,$ it follows that $\dim(X^{\perp_{Q_i}})=n+2$ for $i=1,2.$ Suppose for a contradiction that $\dim (T_{\pX} B) \geq n+1.$ Then $$X^{\perp_{Q_1}} = X^{\perp_{Q_2}} =: H.$$ Since the cone points span the entire $\bbp(U)$, there exists a cone point $[v_\lambda]$ of a singular quadric $Q_\lambda\in\scl$ such that $v_\lambda\notin H.$ Since $Q_\lambda$ descents to a quadratic form on the 2-dimensional vector space $H/X$, there exists a vector $v\in H\backslash X$ such that $Q_\lambda(v)=0.$ Now, $$\tSpan\{X,v,v_\lambda\}\subset U$$ is an $(n+2)$-dimensional isotropic subspace with respect to $Q_\lambda.$ However, since $Q_\lambda$ is a simple quadric cone, its maximal isotropic subspace has dimension $n+1.$\epr
\end{example}

We will prove Theorem \ref{thm:Hsimpletransitive} by proving the following three Propositions.

\begin{proposition}\label{prop:transitive} $\tDiv(C)$ acts transitively on $F\dcup F'.$
\end{proposition}

\begin{proposition}\label{prop:principaldies} The principal divisors act trivially on $F\dcup F'.$ Since $[D_0]$ acts trivially, we now have a transitive action of $H$ on $F\dcup F'.$
\end{proposition}

\begin{proposition}\label{prop:nothingelsedies} If $[D]\in H$ acts trivially, then $[D]=0.$
\end{proposition}

Without loss of generality, we assume that $k$ is algebraically closed. The following two lemmas proved in \cite{Donagi} are crucial in proving these propositions.

\begin{lemma}\label{lem:Don26} (Lemma 2.6 in \cite{Donagi}) Suppose $\pX,\pX'\in F$ intersect at codimension $r$. There exists a unique effective divisor $D$ of degree $r$ such that
$$\pX + D = \pX'\mbox{ if }r\mbox{ is even},\quad\pX + D = -\pX'\mbox{ if }r\mbox{ is odd}.$$
In particular, there exists $\pX''\in F$ intersecting $\pX$ at codimension 1 and $\pX'$ at codimension $r-1.$
\end{lemma}

\begin{lemma}\label{lem:Don32} (Lemma 3.2 in \cite{Donagi}) Suppose $[D]\in\mbox{Pic}(C)$ is effective with $h^0(D)=\dim H^0(\sco_C[D])\geq2$, where $$H^0(\sco_C[D]) = \{f\in k^a(C)|[D] + \tdiv(f)\geq 0\}.$$ Then $[D]-[D_0]$ is also effective.
\end{lemma}

\textbf{Proof of Proposition \ref{prop:transitive}:} It suffices to show the existence of an element $D\in\tDiv(C)$ sending $-\pX$ to $\pX'$ for both $\pX,\pX'\in F.$ First suppose $T_{\pX} B$ is an $n$-plane and there exists a point $e\in C$ such that $-\pX+(e)=\pX.$ (cf. Example \ref{ex:bigTxB}) We claim via induction on the codimension $r$ of the intersection $X'\cap X$ in $X,$ that there is an element $D\in\tDiv(C)$ such that $[D]+(-\pX)=\pX'.$ The base case $r=0$ is when $\pX=\pX'$, in which case $[D]=(e)$ does the job. The case $r=1$ is covered by Example \ref{ex:codim1}. Suppose the claim is true for all $\pX''$ intersecting $\pX$ at codimension $\leq r-1$ and codim$(\pX'\cap \pX)=r.$ Choose any $\pX''\in F$ intersecting $\pX'$ at codimension 1 and $\pX$ at $r-1.$ Denote by $D''\in\tDiv(C)$ the element sending $-\pX$ to $\pX''.$ Consider $$D'=(\pX'+\pX'')-D''+(e).$$
From our definition of the action of $H$ on $F\dcup F',$ we know that $$-D''+\pX=-(D''+(-\pX))=-\pX''.$$
Now $(e)$ sends $-\pX$ to $\pX$, $-D''$ sends $\pX$ to $-\pX'',$ and $(\pX' + \pX'')$ sends $-\pX''$ to $\pX'.$ Therefore the composition $D'$ sends $-\pX$ to $\pX'$ as desired.

Next, let $\pX',\pX''\in F$ be arbitrary. Let $D',D''$ denote the elements in $\tDiv(C)$ sending $-\pX$ to $\pX',\pX''$ respectively. Consider
$$D=D'-(e)+D''.$$
Now $D''$ sends $-\pX''$ to $\pX,$ $-(e)$ sends $\pX$ to $-\pX,$ and $D'$ sends $-\pX$ to $\pX'.$ \epr

\begin{lemma}\label{lem:onefiximpltriv} If $D\in\tDiv(C)$ fixes some $x_0\in F\dcup F',$ then $D$ acts trivially.
\end{lemma}

\bpr This follows immediately from the transitivity of the action.\epr

\begin{lemma}\label{lem:DEleqn} If $D,E$ are effective divisors of degree at most $n$, and $D-E=\tdiv(f)$ is a principal divisor, then $D-E$ acts trivially.
\end{lemma}

\bpr Applying Lemma \ref{lem:Don32} repeatedly to $D$, one obtains an unique effective divisor $D_1$ with $h^0(D_1)=1$ and such that $D$ and $E$ are in the linear system $D_1+\frac{\deg(D) - \deg(D_1)}{2}D_0.$ Since $\deg(D)\leq n,$ $H^0(\sco_C(\frac{\deg(D) - \deg(D_1)}{2}D_0)$ consists of functions pulled back from $\bbp^1.$ Hence $D-E$ is a linear combination of divisors of the form $(P)+(\bar{P})$ which acts trivially on $F\dcup F'$ by construction.\epr

Let $\infty$ denote a Weierstrass point of $C$ defined over $k^a.$

\begin{lemma}\label{lem:Ddescent} Suppose $D = (P_1) + \cdots + (P_r) - r(\infty)\in\tDiv(C)$ with $P_i\neq\infty$ and $r\leq n.$ If $D$ is linearly equivalent to $E = (Q_1) + \cdots + (Q_{r'}) - r'(\infty)$ with $Q_i\neq\infty$ and $r'\leq r,$ then $x + D = x + E$ for all $x\in F\dcup F'.$
\end{lemma}

\bpr Apply Lemma \ref{lem:DEleqn} to the effective divisors $(P_1) + \cdots + (P_r)$ and $(Q_1) + \cdots + (Q_{r'}) + (r-r')(\infty).$\epr

Every divisor class $[D]\in J=\picz(C)$ can be represented by a divisor of the form $(P_1) + \cdots + (P_r) - r(\infty)$ with $r\leq n.$ Lemma \ref{lem:Ddescent} says that two different representations of $[D]$ have the same action on $F\dcup F'.$ Since $\deg(D)$ is even, it sends $F$ to $F$. Therefore we have a morphism of varieties $$\al:J\rightarrow\mbox{Aut}(F).$$ The image of $\al$ lies in a commutative subvariety of $\mbox{Aut}(F).$ Since $J$ is complete and $\al([0])=\text{id},$ rigidity\footnote{Rigidity lemma (\cite{Mumford}, pp.40--41): $X$ complete, $Y,Z$ any variety, $f:X\times Y\rightarrow Z$ a morphism such that $f(X\times\{y\})=\{z\}$ for some $y\in Y$. Then $f$ factors as $X\times Y\xrightarrow{p_2}Y\xrightarrow{g} Z.$} implies that $\al$ is a group homomorphism.

\textbf{Proof of Proposition \ref{prop:principaldies}:} Let $\be:\tDiv^0(C)\rightarrow\mbox{Aut}(F)$ denote the action map. To show the principal divisors act trivially, it suffices to show $\be$ factors through $\al:J\rightarrow\mbox{Aut}(F).$ Both are group homomorphisms, therefore it suffices to check $$\be((c)-(c'))=\al([(c)-(c')])$$ for any $c,c'\in C.$ For any $\pX\in F,$
\begin{eqnarray*}
\al([(c)-(c')])(\pX) &=& \pX + (c)-(\infty) + (\bar{c'})-(\infty) \\
&=& \pX + (c) - (c') \\
&=& \be((c)-(c'))(\pX).\epr
\end{eqnarray*}

Given two elements $x=\pm\pX,x'=\pm\pX'$ of $F\dcup F',$ we define their \textbf{intersection codimension} as the intersection codimension of $\pX,\pX'$ and write $$\mbox{codim}(x,x')=\mbox{codim}(\pX,\pX').$$ In this notation, Lemma \ref{lem:Don26} can be stated as follows:

\begin{lemma}\label{lem:Don26re} Suppose $x,x'\in F$ or $x,x'\in F'.$ Then there exists a unique effective divisor $D$ of degree $r=\mbox{codim}(x,x')$ such that $$x + D = (-1)^rx'.$$
\end{lemma}

\begin{lemma}\label{lem:degcodim} Suppose $D$ is an effective divisor of degree $r\leq n,r\geq 1$, then there exists an $x\in F$ such that $$\mbox{codim}(x,x+D)\equiv r\pmod{2}.$$ There is also an $x\in F'$ satisfying the same condition.
\end{lemma}

\bpr Suppose for a contradiction that for all $x\in F,$
\begin{equation}\label{eq:contra}
\mbox{codim}(x,x+D)\equiv r-1\pmod{2}.
\end{equation}
Consider the following variety
$$\Sigma = \{(x,c_1,\ldots,c_{r-1})|x\in F,c_i\in C, x + D = -x + (c_1) + \cdots + (c_{r-1})\}\subset F\times \mbox{Sym}^{r-1}(C).$$
When $r=1$, $\Sigma=\{x\in F|x + D = -x\}.$ It is clear from the definition that $\Sigma$ is closed. Denote the two projections to $F$ and $\mbox{Sym}^{r-1}(C)$ by $\pi_1,\pi_2$ respectively. For any $x\in F,$ $$\mbox{codim}(x,x+D)=:r'\leq r.$$ By Lemma \ref{lem:Don26re}, there exists an effective divisor $D'$ of degree $r'$ such that $$x + D = (-1)^{r-r'}(x+D').$$ Assumption \eqref{eq:contra} says $r-r'$ is odd for all $x$. Therefore, replacing $D'$ by $D'+(r-1-r')(\infty),$ we see that $\pi_1$ is surjective. Since $\mbox{dim}(F)\geq n$ and $\mbox{dim}(\mbox{Sym}^{r-1}(C))=r-1<n,$ there exists a fiber of $\pi_2$ of positive dimension. In other words, there exists a divisor $\w{D}$ of odd degree such that for infinitely many $x\in F$,
\begin{equation}\label{eq:double}
x + \w{D} = -x.
\end{equation}

Let $D_1$ be a divisor such that $2D_1 - (\infty)$ is linearly equivalent to $\w{D}.$ Since we have shown that the principal divisors act trivially, \eqref{eq:double} implies that for infinitely many $x\in F,$ $$(x+D_1) = -(x+D_1) + (\infty).$$ Hence for infinitely many $\pX\in F,$ $$\pX = \tau(\infty)\pX.$$ However, as we will see in Example \ref{ex:Weier}, there are only $2^{2n}$ such $\pX.$ Contradiction.

The statement for $F'$ follows from the same argument, which is the main reason why we have used $x$ to denote an element of $F$ instead of the usual $\pX.$\epr

\textbf{Proof of Proposition \ref{prop:nothingelsedies}:} Suppose $D = (P_1)+\cdots + (P_r)-r(\infty)$ acts trivially on $F$ with $r\leq n$ minimal and $P_i\neq\infty$.

Suppose first that $r=2r'$ is even. Then for all $\pX\in F, $
\begin{equation}\label{eq:evenr}
\pX + (P_1)+\cdots + (P_{r'}) = \pX + (\bar{P}_{r'+1})+\cdots + (\bar{P}_r).
\end{equation}
By lemma \ref{lem:degcodim}, there exists $\pX_0\in F$ such that $$\mbox{codim}(\pX_0,\pX_0 + (P_1)+\cdots + (P_{r'}))=r''\equiv r'\pmod{2}.$$ Therefore, there exists points $Q_1,\ldots,Q_{r''}\in C$ such that $$\pX_0 + (P_1)+\cdots + (P_{r'}) = \pX_0 + (Q_1)+\cdots + (Q_{r''}).$$ Lemma \ref{lem:onefiximpltriv} says if a divisor fixes one $\pX_0\in F,$ then it acts trivially on $F$. Hence the divisor $$(Q_1)+\cdots + (Q_{r''}) + (P_{r'+1})+\cdots + (P_{r})-(r''+r')(\infty)$$ acts trivially on $F.$ Minimality of $r$ forces $r''=r'.$ That is, $$\mbox{codim}(\pX_0,\pX_0 + (P_1)+\cdots + (P_{r'}))=r'.$$ Lemma \ref{lem:Don26} then implies $$(P_1)+\cdots + (P_{r'}) = (\bar{P}_{r'+1})+\cdots + (\bar{P}_r)$$ as effective divisors of degree $r'$. Therefore $D=0.$

Suppose now $r = 2r'+1$ is odd. Then for all $\pX\in F, $
\begin{equation}\label{eq:oddr}
\pX + (P_1)+\cdots + (P_{r'+1}) = \pX + (\bar{P}_{r'+2})+\cdots + (\bar{P}_r) + (\infty)
\end{equation}
Argue just like the even case, we see that minimality of $r$ implies that for some $\pX_0\in F,$ $$\mbox{codim}(\pX_0,\pX_0 + (P_1)+\cdots + (P_{r'+1}))=r'+1.$$ Then Lemma \ref{lem:Don26} implies $$(P_1)+\cdots + (P_{r'+1}) = (\bar{P}_{r'+2})+\cdots + (\bar{P}_r)+(\infty)$$ as effective divisors of degree $r'+1$. Therefore $D=0.$\epr

We have completed the proofs of Propositions \ref{prop:transitive}, \ref{prop:principaldies}, and \ref{prop:nothingelsedies}. Before moving on to state the main theorem, we describe a stronger form of Lemma \ref{lem:degcodim} for completeness.

Lemma \ref{lem:Don32} implies that if $(P_1)+\cdots+(P_r)-r(\infty)$ and $(Q_1)+\cdots+(Q_r)-r(\infty)$, with $r\leq n$, are two distinct divisors representing the same divisor class $[D]\in J,$ then $[D]$ can also be represented by a divisor of the form $(R_1)+\cdots+(R_{r-2})-(r-2)(\infty).$ Therefore if $r$ is minimal among all such representations of $[D]$, there is a unique effective divisor $D'$ of degree $r$ such that $$[D'-r(\infty)]=[D].$$ We call $D'$ the \textbf{$\infty-$minimal form} of $[D].$

\begin{cor} Let $D'$ be the $\infty-$minimal form of a nonzero divisor class $[D].$ Then there exists an $x\in F$ such that $$\mbox{codim}(x,x+D')=\deg(D').$$ There is also an $x\in F'$ satisfying the same condition.
\end{cor}

\bpr Let $r$ denote the degree of $D'.$ Lemma \ref{lem:degcodim} allows us to pick an $x\in F$ such that $$\mbox{codim}(x,x+D')=:r'\equiv r\pmod{2}.$$ By Lemma \ref{lem:Don26re}, there exists an effective divisor $D''$ of degree $r'$ such that $x + D' = x + D''.$ Hence $D'-D''$ fixes $x$ and by Lemma \ref{lem:onefiximpltriv}, $D'-D''$ acts trivially on $F$. By Proposition \ref{prop:nothingelsedies}, $D'$ is linearly equivalent to $D''$. Since $D'$ is the $\infty-$minimal form of $[D]$, we see that $r'=r.$

The statement for $F'$ follows from the same argument.\epr

We now state our theorem in its completion. 

\begin{theorem}\label{thm:maintheorem}
There is a unique commutative algebraic group structure $+_G$ defined over $k$ on the disconnected variety $$G = \picz(C) \dcup F \dcup \pico(C) \dcup F'$$ such that,
\begEnu
\item $G^0 = \picz(C)$ with component group $G/G^0 \simeq \bbz/4,$
\item $F'$ is isomorphic to $F$ as varieties via the inversion map $-1_G$,
\item the group law extends that on $H = \Pic{C}/D_0 \simeq \picz(C) \dcup \pico(C)$ where $D_0$ is the hyperelliptic class,
\item the group law defines a simply-transitive action of $H$ on $F\dcup F'$ extending the following action of $C:$
$$\pX+(c)=-\tau(c)\pX,\qquad -\pX+(c)=\tau(\bar{c})\pX,$$
with respect to which $x+_Gx',$ for $x,x'\in F\dcup F'$, is the unique divisor class sending $-x$ to $x'$.
\endEnu
\end{theorem}

\bpr The only thing left to check is associativity, which amounts to the following four:
\begin{eqnarray*}
[D_1]+_G([D_2]+_G[D_3])&=&([D_1]+_G[D_2])+_G[D_3]\\
x+_G([D_2]+_G[D_3])&=&(x+_G[D_2])+_G[D_3]\\
x+_G(x'+_G[D_3])&=&(x+_Gx')+_G[D_3]\\
x+_G(x'+_Gx'')&=&(x+_Gx')+_Gx'',
\end{eqnarray*}
for $[D_1],[D_2],[D_3]\in H$ and $x,x',x''\in F\dcup F'.$

The first one is associativity of the group law on $H$. The second
follows from the definition of the action of $H$. The third follows
as both sides send $-x$ to $x'+[D_3].$ For the fourth one, denote
the two sides by $x_L$ and $x_R$ and add $x'$ to both sides. The
third associativity tells us $x'+_Gx_L=(x'+_Gx)+_G(x'+_Gx'')$ and likewise, $x_R+_Gx'=(x+_Gx')+_G(x''+_Gx')$. Commutativity of $+_G$ implies these two elements of $\picz(C)$ are equal. Therefore $x_L=x_R$ is the image of $-x'.$ \epr

\begin{cor}
The class $[F]\in H^1(k,J)$ is 4-torsion, twice of which is $[\pico(C)].$ One can lift $[F]$ to a torsor of $J[4]$ by taking $$F[4]:=\{\pX\in F| \pX+_G\pX+_G\pX+_G\pX=0\}.$$
\end{cor}

\bpr With our convention of Galois cohomology, we need to show $F(k^s)$ is nonempty. Let $P$ be a Weierstrass point defined over $k^s.$ There are precisely $2^{2n}$ elements of $F(k^a)$ satisfying $\pX +_G \pX = (P).$ As we will see in Example \ref{ex:Weier}, they correspond to $(n-1)$-planes contained in the base locus of a generic pencil in $\bbp^{2n}.$ Theorem \ref{thm:simpletransitivesmoothclose} says they are in fact all defined over $k^s$.\epr

When $\pico(C)(k)\neq\emptyset,$ $[F]$ is 2-torsion. For each $[D_1]\in\pico(C)(k),$ we obtain a lift of $F$ to a torsor of $J[2]$ by taking
$$F[2]_{[D_1]}=\{\pX\in F|\pX +_G \pX = [D_1]\}.$$

\subsubsection{Example: rational Weierstrass point}
\begin{example}\label{ex:Weier}
Suppose $C$ has a rational Weierstrass point, or equivalently, $\scl$ has a rational singular quadric. By moving this point to $\infty,$ we assume that $Q_1$ is singular with cone point $[v_\infty].$ Let $H=v_\infty^{\perp_{Q_2}}$ be the hyperplane in $U$ orthogonal to $v_\infty$ with respect to $Q_2$. Then $\tau(\infty)$ is induced by the linear map on $U$ that fixes $H$ and sends $v_\infty$ to $-v_\infty.$ Hence,
\begin{equation}\label{eq:F2weier}
F[2]_\infty = \{\pX\in F|\pX\subset B\cap \bbp H\}.
\end{equation}

Notice when restrict to the $2n+1$ dimensional vector space $H$, $Q_1$ and $Q_2$ span a generic pencil $\scl_H$. Moreover, $Q_{1|H}$ is non-degenerate. Let $T$ be the self-adjoint operator on $H$ associated to the pencil $\scl_H$ as defined in \eqref{eq:defT}. Then the right hand side of \eqref{eq:F2weier} is precisely $W_T$ as defined in the odd dimension case. Now $J[2]$ acts on $F[2]_\infty$ via the action of $J$ and on $W_T$ via the identification $J[2]\simeq \tStab_{\PO(H,Q_{1|H})}(T).$ The following Proposition says that these two actions coincide.

\begin{proposition}\label{prop:FT2WT} $F[2]_\infty = W_T$ as $J[2]$-torsors.
\end{proposition}

\vspace{5pt}\bpr It suffices to show for any $(P)-(\infty)\in J[2](k^s)$ with $P$ a Weierstrass point, the two actions are the same. Let $\al$ denote the root of $f(x)$ corresponding to $P$, and set $h(x) = f(x)/(x-\al).$ On $W_T(k^s),$ by Remark \ref{rmk:J2inPO}, the action of $(P)-(\infty)$ is induced by the following map on $H\otimes k^s:$
$$x\mapsto x - 2\frac{h(T)}{h(\al)}x.$$

We now compute the action of $(P)-(\infty)$ on $F[2]_\infty(k^s).$ The singular quadric corresponding to $P$ is $\al Q_1 - Q_2.$ Let $w_P\in H\otimes k^s$ be an eigenvector of $T$ with eigenvalue $\al.$ The cone point of $\al Q_1 - Q_2$ is $[(w_P,0)].$ Here we have decomposed $U\otimes k^s$ as $H\oplus U_{\infty,T}$ where $U_{\infty,T}$ is the kernel of the degenerate quadric $Q_1$. Let $b_1,b_2$ denote bilinear forms associated to $Q_1,Q_2$ and also to their restriction to $H$. From the definition of $\tau$ earlier, cf. \eqref{eq:defoftau}, we see that the action of $(P)-(\infty)$ is induced by the following map on $U\otimes k^s:$
$$x\mapsto x - 2\frac{b_2(x,(w_P,0))}{b_2((w_P,0),(w_P,0))}(w_P,0).$$
If we view each $\pX\in F[2]_\infty(k^s)$ as sitting inside $\bbp(H),$ then the action of $(P)-(\infty)$ is induced by the following map on $H\otimes k^s:$ $$x\mapsto x-2\frac{b_1(x,w_P)}{b_1(w_P,w_P)}w_P.$$

To prove the lemma, it remains to show for any $x\in H\otimes k^s,$ $$\frac{h(T)}{h(\al)}x = \frac{b_1(x,w_P)}{b_1(w_P,w_P)}w_P.$$ Since both sides are killed by $T-\al,$ and since $T$ has 1-dimensional eigenspaces, they are both scalar multiples of $w_P.$ Now
$$b_1(\frac{h(T)}{h(\al)}x,w_P)=b_1(x,\frac{h(T)}{h(\al)}w_P)=b_1(x,w_P)=b_1(\frac{b_1(x,w_P)}{b_1(w_P,w_P)}w_P,w_P).$$ Therefore they are the same scalar multiple of $w_P.$\epr

\begin{remark}\label{rmk:J2inPOalt}
Equation \eqref{eq:F2weier} offers another view point for the canonical identification of $J[2]$ with the stabilizer of a self-adjoint operator, namely they share a common principal homogeneous space. Fix any $k$-rational $T$, then $J[2]$ acts on $F[2]_\infty$ simply-transitively and $\tStab(T)$ acts on $W_T$ simply-transitively. It is clear from the definitions that these two actions commute. Fix some $X_0\in F[2]_\infty$, one can define the map
$$\iota:J[2]\rightarrow \tStab(T)$$
by taking $\iota([D])$, for any $[D]\in J[2]$, to be the unique element of $\tStab(T)$ sending $X_0$ to $X_0 + [D].$ Commutativity of the two actions and commutativity of $J[2]$ show that this map is independent on the choice of $X_0.$ Proposition \ref{prop:FT2WT} then implies that $\iota$ is given by the map we defined in Remark \ref{rmk:J2inPO}.
\end{remark}
\end{example}

\subsubsection{Example: rational non-Weierstrass point}
\begin{example}\label{ex:nonWeier}
Suppose $C$ has a rational non-Weierstrass point, or equivalently, $\scl$ has a rational quadric with discriminant 1. By moving the point to infinity, we assume that $Q_1$ has discriminant 1. Its two rulings are therefore defined over $k$. Let $Y_0$ denote one of the rulings and let $\infty\in C(k)$ denote the point corresponding to the quadric $Q_1$ and the ruling $Y_0$. Denote by $\infty'$ the conjugate of $\infty$ under the hyperelliptic involution. Let $T$ denote the self-adjoint operator on $U$ associated to the pencil $\scl$ as defined in \eqref{eq:defT}.

\begin{proposition}\label{prop:reductionEven} $$F[2]_\infty = \{\pX\in F|\pX = \tau(\infty)\pX\} = \{\pX\simeq\bbp^{n-1} | \tSpan\{\pX,\bbp(TX)\}\sim Y_0\}.$$
The latter condition means $\tSpan\{\pX,\bbp(TX)\}\simeq\bbp^n$ is contained in $Q_1$ in the ruling $Y_0$.
\end{proposition}

\bpr Suppose $\pX\simeq\bbp^{n-1}$ with $\tSpan\{\pX,\bbp(TX)\}\sim Y_0.$

1). Since $TX\subset X^{\perp_{Q_1}},$ we see $X\subset X^{\perp_{Q_2}}$ and hence $\pX\in F.$

2). Since $\tSpan\{\pX,\bbp(TX)\}\supset \pX$ is an $n$-plane contained in $Q_1$ in the same ruling as $\pY_0,$ we see $\tau(\infty)\pX$ is the residual intersection of $\tSpan\{\pX,\bbp(TX)\}$ with $Q_2$.

3). $\tSpan\{\pX,\bbp(TX)\}$ intersects $Q_2$ tangentially at $\pX$ because $$TX\subset TX^{\perp_{Q_1}}\impl TX \subset X^{\perp_{Q_2}}\impl T_{\pX}Q_2\supset\tSpan\{\pX,\bbp(TX)\}.$$ Therefore $\pX\in F[2]_\infty.$

Conversely, suppose $\pX\in F[2]_\infty.$ Suppose $\tSpan\{\pX,[p]\}\supset\pX$ is the $n$-plane contained in $Q_1$ in the same ruling as $\pY_0$, for some $p\in U\otimes k^a.$ Since $\tau(\infty)\pX = \pX,$ we see $$b_1(x,p) = b_1(x,Tp) = b_1(p,p) = 0, \forall x\in X.$$

1). Since $\tSpan\{\pX,[p]\}$ does not lie in the base locus, $b_1(p,Tp)=Q_2(p)\neq 0.$

2). Since $\tSpan\{X,p\}\subset p^{\perp_{Q_1}},$ we have $Tp\notin\tSpan\{X,p\}$ but $Tp\in X^{\perp_{Q_1}}$. Hence $$X^{\perp_{Q_1}} = \tSpan\{X,p,Tp\}.$$

3). Since $TX\subset p^{\perp_{Q_1}}\cap X^{\perp_{Q_1}},$ we have $TX\subset\tSpan\{X,p\}.$

4). If $TX\subset X,$ then $X^{\perp_{Q_1}}=X^{\perp_{Q_2}}$ which implies that $T_{\pX}(Q_1\cap Q_2)\simeq \bbp^{n+1}$. This contradicts Lemma \ref{lem:dimTxB}. Therefore $$X\subsetneq TX\subset\tSpan\{X,p\},\quad \mbox{i.e.}\quad \tSpan\{\pX,\bbp(TX)\}=\tSpan\{\pX,[p]\}\sim Y_0.\epr$$

In parallel to the odd dimension case, Proposition \ref{prop:reductionEven} then suggests fixing the monic polynomial $f(x)$ of degree $2n+2$, the quadratic form $Q_1$ of discriminant 1, and considering the following $k$-schemes,
\begin{eqnarray*}
V_f &=& \{T:U\rightarrow U | T^* = T, \mbox{ characteristic polynomial of }T\mbox{ is }f\}\\
W_f &=& \{(T,X)\in V_f\times\mbox{Gr}(n,U)| \tSpan\{X,TX\}\sim Y_0\}
\end{eqnarray*}
Here $\tSpan\{X,TX\}\sim Y_0$ means that $\tSpan\{X,TX\}$ is a (linear) $(n+1)$-plane isotropic with respect to $Q_1$ lying in the ruling $Y_0$. Let $W_T$ denote the fiber of $W_f$ above $T$, then Proposition \ref{prop:reductionEven} says $F[2]_\infty = W_T.$ The group $\PSO(U,Q_1)$ preserves the rulings, hence it acts on $W_f$ via
$$g.(T,X) = (gTg^{-1},gX).$$

\begin{proposition}\label{prop:PSOtransitive}
The action of $\PSO(U,Q_1)$ on $V_f$ has a unique geometric orbit. For any $T\in V_f(k')$ defined over some field $k'$ over $k$, its stabilizer scheme $\tStab(T)$ is isomorphic to $(\mbox{Res}_{L'/k'}\mu_2)_{N=1}/\mu_2\simeq J_{k'}[2]$ as group schemes over $k'$ where $L' = k'[x]/f(x).$
\end{proposition}

\bpr For any $T$ in $V_f(k'),$ since $T$ is regular semi-simple, its stabilizer scheme in $\GL(U_{k'})$ is a maximal torus. It contains and hence equals to the maximal torus $\mbox{Res}_{L'/k'}\bbg_m.$ For any $k'$-algebra $K$,
\begin{eqnarray*}
\tStab_{O(U_{k'},Q_1)}(T)(K) &=& \{g\in(K[T]/f(T))^\times|g^*g = 1\}\\
&=& \{g\in(K[T]/f(T))^\times|g^2 = 1\}.
\end{eqnarray*}
Hence
\begin{eqnarray*}
\tStab_{O(U_{k'},Q_1)}(T) &\simeq& \mbox{Res}_{L'/k'}\mu_2,\\
\tStab_{SO(U_{k'},Q_1)}(T) &\simeq& (\mbox{Res}_{L'/k'}\mu_2)_{N=1},\\
\tStab_{PSO(U_{k'},Q_1)}(T) &\simeq& (\mbox{Res}_{L'/k'}\mu_2)_{N=1}/\mu_2.
\end{eqnarray*}

Suppose $T_1,T\in V_f(k^s)$. There exists $g\in GL(U)(k^s)$ such that $T_1 = gTg^{-1}.$ Since $T_1$ and $T$ are both self-adjoint, $g^*g$ centralizes $T$ and hence lies in $(k^s[T]/f(T))^\times$ which is a product of $k^{s\times}$ since $f$ splits. Since the characteristic of $k$ is not 2, there exists $h\in (k^s[T]/f(T))^\times$ such that $g^*g=h^2.$ Then $gh^{-1}$ is an element of $O(U,Q_1)(k^s)$ conjugating $T$ to $T_1.$ Multiplying the $h$ by $(-1,1,\ldots,1)\in(k^s[T]/f(T))^\times$ if necessary, we may assume $gh^{-1}\in \SO(U,Q_1)(k^s).$ Its image in $\PSO(U,Q_1)(k^s)$ does the job.\epr

For general $Q$, there might not be a self-adjoint operator defined over $k$ with the prescribed characteristic polynomial. For example over $\bbr$, operators self-adjoint with respect to the positive definite form have real eigenvalues.

\begin{lemma}\label{lem:PSOsplit} If $Q_1$ is split, then both $V_f(k)$ and $W_f(k)$ are nonempty. Furthermore, there exists $(T_0,X_0)\in W_f(k)$ with trivial stabilizer in $\PSO(U,Q_1)(k^a).$
\end{lemma}

\bpr Consider the $2n+2$ dimensional \'{e}tale $k$-algebra $L = k[x]/f(x)=k[\be].$ On $L$ there is the following bilinear form
$$<\lambda,\mu> = \tTr(\lambda\mu/f'(\be))=\mbox{ coefficient of }\be^{2n+1}\mbox{ in }\lambda\mu.$$
This form defines a split quadratic form since $Y=\tSpan_k\{1,\be,\ldots,\be^n\}$ is a rational maximal isotropic subspace. Hence there exists an isometry from $(L,<,>)$ to $(U,Q_1)$ defined over $k$. Denote by $T_0$ the image of the multiplication by $\be$ operator, and by $X_0$ the image of $X=\tSpan_k\{1,\be,\ldots,\be^{n-1}\}$. Since $(\cdot\be,X)$ has trivial stabilizer in $\PSO(L,<,>)(k^a),$ its image $(T_0,X_0)$ has trivial stabilizer in $\PSO(U,Q_1)(k^a).$\epr

\begin{theorem}\label{thm:PSOsimpletransitivesmoothclose} Suppose $k$ is algebraically closed or separably closed, then $\PSO(U,Q_1)(k)$ acts simply-transitively on $W_f(k).$
\end{theorem}

\bpr Suppose $k$ is separably closed. Proposition \ref{prop:PSOtransitive} shows it suffices to prove that for the $T_0\in V_f(k)$ obtained in the above lemma, $\tStab(T_0)(k)$ acts simply-transitively on $W_{T_0}(k).$ Since $(T_0,X_0)$ has trivial stabilizer, it suffices to show they have the same size. As a consequence of Proposition \ref{prop:reductionEven}, for any $k$, $W_{T}(k^a)\simeq F[2]_\infty(k^a)\simeq J[2](k^a)$ has $2^{2n}$ elements for any $T$. Hence we are done because,
$$2^{2n} = |(\Rlk\mu_2/\mu_2)_{N=1}(k)| = |\tStab(T_0)(k)| \leq |W_{T_0}(k)| \leq |W_{T_0}(k^a)|=2^{2n}.\epr$$

\begin{cor}\label{thm:PSOsimpletransitivesmooth} $\PSO(U,Q_1)(k')$ acts simply-transitively on $W_f(k')$ for any field $k'$ over $k$.
\end{cor}

\bpr It suffices to prove transitivity. Suppose $(T_1,X_1),(T_2,X_2)\in W_f(k'),$ let $g\in \PSO(U,Q_1)(k'^s)$ be the unique element sending $(T_1,X_1)$ to $(T_2,X_2).$ Then for any $\sigma\in \Gal(k'^s/k'),$ $\gl{g}{\sigma}$ also sends $(T_1,X_1)$ to $(T_2,X_2).$ Hence $g = \gl{g}{\sigma}\in \PSO(U,Q_1)(k').$\epr

\begin{remark}\label{rmk:J2inPSO}
One can write down an explicit formula for the identification,
\begin{equation}\label{eq:PSOStabJ2}
\tStab(T)\simeq (\Rlk\mu_2)_{N=1}/\mu_2\simeq J[2],
\end{equation}
The method is the same as the odd case in Remark \ref{rmk:J2inPO}. Denote the roots of $f(x)$ over $k^s$ by $\al_1,\ldots,\al_{2n+2},$ and for each $i$, define $h_i(x) = f(x)/(x-\al_i).$ Then on the level of $k^s$-points, \eqref{eq:PSOStabJ2} is given by sending $$\sum n_i(\al_i) - \frac{\sum n_i}{2}((\infty)+(\infty')),\quad\sum n_i\mbox{ even},$$ to the image in $\PSO_{2n+2}(k^s)$ of \begin{equation}\label{eq:PSOformula}
\prod\left(1-2\frac{h_i(T)}{h_i(\al_i)}\right)^{n_i} = 1 - 2\sum n_i\frac{h_i(T)}{h_i(\al_i)}.
\end{equation}
Note as a polynomial of degree at most $2n+1,$ $\sum_{i=1}^{2n+2} h_i(x)/h_i(\al_i)$ takes the value 1 when $x=\al_1,\ldots,\al_{2n+2},$ hence it must be the constant polynomial 1. Thus,
$$\prod_{i=1}^{2n+2}\left(1-2\frac{h_i(T)}{h_i(\al_i)}\right)=-1=1\mbox{ in }\PSO_{2n+2}.$$
We will see in Proposition \ref{prop:FT2WT} and Proposition \ref{prop:PSOFT2WT} that $1-2\frac{h_i(T)}{h_i(\al_i)}$ is a reflection, hence has determinant $-1$. The assumption that $\sum n_i$ is even ensures that the product in \eqref{eq:PSOformula} lies in $\SO.$
\end{remark}

Just as in the odd case, $J[2]$ acts on $F[2]_\infty$ via the action of $J$ and on $W_T$ via the identification $J[2]\simeq \tStab_{\PSO(U,Q_1)}(T).$ The following Proposition says that these two actions coincide.

\begin{proposition}\label{prop:PSOFT2WT} $F[2]_\infty = W_T$ as $J[2]$-torsors.
\end{proposition}

\vspace{5pt}\bpr It suffices to show for any $(P_1)-(P_2)\in J[2]$ with $P_1,P_2$ any two Weierstrass points, the two actions are the same. Let $\al_i$ denote the root of $f(x)$ corresponding to $P_i$, and set $h_i(x) = f(x)/(x-\al_i).$ On $W_T(k^a),$ by Remark \ref{rmk:J2inPSO}, the action of $(P_1)-(P_2)$ is induced by the following map on $U\otimes k^s:$ $$x\mapsto x - 2\frac{h_1(T)}{h_1(\al_1)}x - 2\frac{h_2(T)}{h_2(\al_2)}x$$ on

For $i=1,2$, let $w_i\in U\otimes k^s$ be an eigenvector of $T$ with eigenvalue $\al_i.$ The cone point of the singular quadric corresponding to $P_i$ is therefore $[w_i]$. Let $b$ denote the bilinear form associated to $Q$. Then on $F(k^s),$ the action of $\tau(P_i)$ is induced by the following map on $U\otimes k^a:$
$$\mbox{refl}_{P_i}:x\mapsto x - 2\frac{b_1(x,v_i)}{b_1(v_i,v_i)}v_i.$$
Composing two such reflections, we see that the action of $\tau(P_1)\tau(P_2)$ is induced by the following map on $U\otimes k^a:$
$$x\mapsto x - 2\frac{b_1(x,w_1)}{b_1(w_1,w_1)}w_1 - 2\frac{b_1(x,w_2)}{b_1(w_2,w_2)}w_2 + \frac{4b_1(x,w_1)b_1(w_1,w_2)}{b_1(w_1,w_1)b_1(w_2,w_2)}w_2.$$
Since self-adjoint operators have pairwise orthogonal eigenspaces, the last term is 0. Also as in the proof of Proposition \ref{prop:FT2WT}, $$\frac{h_i(T)}{h_i(\al_i)}x = \frac{b_1(x,w_i)}{b_1(w_i,w_i)}w_i.$$
Therefore the two actions are equal.\epr

\begin{remark}\label{rmk:J2inPSOalt}
In parallel to the odd case, the equality $F[2]_\infty=W_T$ as $\Gal(k^s/k)$-sets provides a different view point on the identification of $J[2]$ with $\tStab(T),$ as they share a common principal homogeneous space. Proposition \ref{prop:PSOFT2WT} implies that this new identification coincides with the formula given by Remark \ref{rmk:J2inPSO}.
\end{remark}
\end{example}


%% file: PenRegular.tex
For the rest of the paper, we focus on regular pencils. Let $\Omega_1$ denote the discriminant hypersurface in $\bbp(H^0(\sco_{\bbp^{N-1}}(2)))$ parameterizing singular quadrics, and let $\Omega_2$ denote the subvariety parameterizing quadrics with higher degeneracy degree. Recall a pencil $\scl$ is generic if and only if $\scl$ is a generic line, and hence it intersects $\Omega_1$ transversely at $N$ points and misses $\Omega_2.$ A pencil is \textbf{regular} if it misses $\Omega_2$ but is allowed to intersect $\Omega_1$ tangentially.

In the $(Q_1,T)$-setup, where $Q_1$ is non-degenerate and $T$ is self-adjoint with respect to $Q_1$, regularity of the pencil is equivalent to regularity of $T$. An operator $T$ is regular if and only if its characteristic polynomial coincide with its minimal polynomial if and only if all its eigenspaces are 1-dimensional. Let $f(x) = \det(xI - T)$ denote the minimal polynomial of $T$ and as before we assume that $f(x)$ splits completely over $k^s$.

The following reduction step is key to study the variety of maximal linear spaces contained in the base locus over $k^s$. Suppose temporarily $k$ is separably closed, let $U$ denote the underlying $N$-dimensional $k$-vector space and let $v\in U$ denote an eigenvector of $T$ whose eigenvalue $\al$ has multiplicity at least 2. Since $v$ is an eigenvector and $T$ is self-adjoint, $T$ descends to a linear operator $\bar{T}$ on $\bar{U}=v^\perp/v$ where $\perp$ is taken with respect to $Q_1$. The quadratic form $Q_1$ descends to a non-degenerate quadratic form $\bar{Q}_1$ on $\bar{U}$ with respect to which $\bar{T}$ is regular self-adjoint with minimal polynomial $f(x)/(x-\al)^2.$ Suppose $N=2n+1$ or $2n+2$ and let $X$ be an $n$-plane in $U$ of interest such that $X\subset X^\perp,TX\subset X^\perp.$ Define $\bar{X}$ to be the image of $X\cap v^\perp$ in $v^\perp/v.$ Then $\bar{X}\subset \bar{X}^\perp,\bar{T}\bar{X}\subset \bar{X}^\perp.$ As we will see in what follows, either $v\in X$ or $X\nsubseteq v^\perp$, hence $\dim\bar{X} = n-1.$ The strategy will be to apply this reduction repeatedly until $T$ becomes regular semi-simple where one can use the result in the previous section regarding generic pencils.

Factors $f(x)$ as $f(x) = \prod_{i=1}^{r+1}(x - \al_i)^{m_i}$ over $k^s$. Let $U_{i,T}$ denote the generalized eigenspace over $k^s$ of $T$ with eigenvalue $\al_i,$ and let $v_i$ be an eigenvector of $T$ with eigenvalue $\al_i.$ Regularity says $v_i$ is unique up to scalars. The singular locus of $B$ consists precisely of $[v_i]$ with $m_i\geq 2.$ For any linear $n$ plane $X$ such that $X\subset X^\perp,TX\subset X^\perp$, where $\perp$ is taken with respect to $Q_1$, and for each $i=1,\ldots,r+1,$ we define $\dim_{i,T}(X)$ to be the dimension of the maximal $T$-stable subspace of $(X\otimes k^s)\cap U_{i,T}.$ Since each $U_{i,T}$ is $m_i$ dimensional and $Q_1$ restricts to a non-degenerate quadratic form on $U_{i,T}$, we have $$\dim_{i,T}(X) \leq m_i/2.$$
For any sequence of integers $d_1,\ldots,d_{r+1}$ such that $0\leq d_i\leq m_i/2,$ we define for any field $k'$ containing $k$,
$$L^{f,T}_{\{d_1,\ldots,d_{r+1}\}}(k')=\{X\simeq (k')^n|X\subset X^\perp,TX\subset X^\perp, \dim_{i,T}(X)= d_i\}.$$
The superscript $f$ is unnecessary, but it serves in making the reduction step clearer. Note $L^{f,T}_{\{0,0,\ldots,0\}}(k')$ is the set of $k'$-points of the variety of projective $n-1$ planes contained in the smooth part of $B$.

\subsection{Odd dimension}

Suppose $N=2n+1$ is odd. For ease of notation, we write $Q$ for the non-degenerate quadratic form $Q_1$. By multiplying $Q$ by a constant, we also assume that $Q$ has discriminant 1. Fixing the minimal polynomial $f(x)$ of degree $2n+1$, let $C$ be the hyperelliptic curve defined by affine equation $y^2 = f(x).$ We define the $k$-scheme,
$$V_f = \{T:U\rightarrow U|T\mbox{ is self-adjoint and regular with characteristic polynomial }f(x)\}.$$
Note here regularity means there is no linear relations between $1,T,\ldots,T^{2n}.$ For every field $k'$ containing $k$, and every $T\in V_f(k'),$ let $W_T(k')$ denote the set of (linear) $n$-dimensional $k'$-subspaces $X$ of $U\otimes k'$ such that $X\subset X^\perp, TX\subset X^\perp.$ As before, we define
$$W_f(k')=\{(T,X)|T\in V_f(k'),X\in W_T(k')\}.$$ There is a Galois invariant action of $\PO(U,Q)=O(U,Q)/(\pm1)$ on $W_f:$ $$g.(T,X)=(gTg^{-1},gX).$$
For any sequence of integers $d_1,\ldots,d_{r+1}$ such that $0\leq d_i\leq m_i/2,$ we define
$$W^f_{\{d_1,\ldots,d_{r+1}\}}(k') = \{(T,X)|T\in V_f(k'),X\in L^{f,T}_{\{d_1,\ldots,d_{r+1}\}}(k')\}.$$

The main theorem we are heading towards is the following:

\begin{theorem}\label{thm:singular}
$|L^{f,T}_{\{d_1,\ldots,d_{r+1}\}}(k^s)|=2^r/2^a,$ where $a$ is the number of $d_i$'s equal to $m_i/2.$
\end{theorem}

The action of $\PO(U,Q)$ preserves the decomposition of $U\otimes k^s$ into generalized eigenspaces, in the sense that
$$U_{i, gTg^{-1}}=gU_{i,T},\quad\quad\forall T\in V_f(K),\forall g\in \PO(U,Q)(K),\forall i=1,\ldots,r+1.$$ Therefore one obtains a Galois equivariant action of $\PO(U,Q)$ on $W^f_{\{d_1,\ldots,d_{r+1}\}}.$

\begin{theorem}\label{thm:simpletransitivesingular}
$\PO(U,Q)(k^s)$ acts on $W^f_{\{d_1,\ldots,d_{r+1}\}}(k^s)$ simply-transitively if $a=0$ and transitively if $a>0.$
\end{theorem}

\begin{cor}\label{thm:simpletransitivesmoothregular} For any field $k'$ over $k$, $\PO(U,Q)(k')$ acts simply-transitively on $W^f_{\{0,\ldots,0\}}(k')$.
\end{cor}

\bpr Same descent argument as in the proof of Corollary \ref{thm:PSOsimpletransitivesmooth}.\epr

We begin by studying the conjugation action of $\PO(U,Q)$ on $V_f.$

\begin{proposition}\label{prop:POtransitiveregular}
The action of $\PO(U,Q)$ on $V_f$ has a unique geometric orbit. For any $T\in V_f(k')$ defined over some field $k'$ over $k$, its stabilizer scheme $\tStab(T)$ is isomorphic to $\mbox{Res}_{L'/k'}\mu_2/\mu_2\simeq J_{k'}[2]$ as group schemes over $k'$ where $L' = k'[x]/f(x).$ In particular, $\tStab_{PO(U,Q)}(T)(k^s)$ is an elementary abelian 2-group of order $2^r,$ where $r+1$ is the number of distinct roots $f(x)$ over $k^s$.
\end{proposition}

\bpr The first statement follows in the same way as in the proof of Proposition \ref{prop:PSOtransitive} except now $k^s[x]/f(x)$ is a product of algebras of the form $k^s[x]/x^{m_i}.$ Every unit in $k^s[x]/x^{m_i}$ is a square as $\mbox{char}(k)\neq2.$

The second statement follows from the structure theory of finitely generated modules over Principal Ideal Domains. One can view $U\otimes k'$ as a module over $k'[x]$ with $x$ acting via the operator $T$. The elements in $\mbox{GL}(U)(k)$ commuting with $T$ are precisely the automorphisms of $U$ as $k'[x]$-modules. Since $T$ is regular, the structure theory of finitely generated modules over PID says that $U\otimes k'$ is isomorphic to $k'[x]/f(x)$ as a $k'[x]$-module. As a module of $k'[x]$ generated by the element 1, the automorphisms of $U$ are precisely multiplication by elements in $(k'[x]/f(x))^\times.$ Then as in Proposition \ref{prop:PSOtransitive},
\begin{eqnarray*}
\tStab_{O(U,Q)}(T)(k') &=& \{g(T)|g\in k'[x], g(T)^*g(T) = 1\}\\
&=&\mu_2(k'[T]^\times)\\
\tStab_{PO(U,Q)}(T)(k') &=&\mu_2(k'[T]^\times)/(\pm1).
\end{eqnarray*}

For the last statement, from the factorization of $f(x)$, we know $$k^s[x]/f(x) \simeq \prod_{i=1}^{r+1}k^s[x]/(x-\al_i)^{m_i}.$$ Therefore, $\tStab_{O(U,Q)}(T)(k^s)\simeq (\bbz/2\bbz)^{r+1}$ is an elementary abelian 2-group of order $2^{r+1}$. Moding out the diagonally embedded $\bbz/2\bbz$ gives $\tStab_{PO(U,Q)}(T)(k^s)$.\epr

\begin{remark}\label{rmk:J2inPOregular}
Just as in Remark \ref{rmk:J2inPO}, we can give a more explicit description for the stabilizer as polynomials in $T$. For each $i=1,\ldots,r+1,$ define $h_i^T(x) = f(x)/(x-\al_i)^{m_i}.$ Then
\begin{eqnarray*}
\mu_2(K[T]^\times) &=& \left\{\prod_{i\in I}\left(1-2\frac{h_i^T(T)}{h_i^T(\al_i)}\right)
\right\}_{I\subset\{1,\ldots,r+1\}}\\
&=&\left\{1-2\sum_{i\in I}\frac{h_i^T(T)}{h_i^T(\al_i)}\right\}_{I\subset\{1,\ldots,r+1\}}.
\end{eqnarray*}
For any $I\subset\{1,\ldots,r+1\}$ and any $j\notin I,$ since $(x-\al_j)^{m_j}$ divides $h_i(x)$ in $K[x]$ and $(T-\al_j)^{m_j}$ kills all the generalized eigenspaces $U_{j,T}$, $$1-2\sum_{i\in I}\frac{h_i^T(T)}{h_i^T(\al_i)}$$ acts trivially on $U_{j,T}$.
\end{remark}

\begin{cor}
For any $T,T'\in V_f(k^s)$, there exists a bijection $$L^{f,T}_{\{d_1,\ldots,d_{r+1}\}}(k^s)\longleftrightarrow L^{f,T'}_{\{d_1,\ldots,d_{r+1}\}}(k^s).$$ 
\end{cor}

\bpr Suppose $g\in \PO(U,Q)(k^s)$ conjugates $T$ to $T'$, then the left action by $g$ on $\mbox{Gr}(n,U)$ gives the desired bijection.\epr

For any $T\in V_f(k^s),$ its stabilizer $J_T$ in $\PO(U,Q)(k^s)$ acts on $L^{f,T}_{\{d_1,\ldots,d_{r+1}\}}(k^s)$. We rephrase the main theorems as follows.

\begin{theorem}\label{thm:regularcool}
For any $X \in L^{f,T}_{\{d_1,\ldots,d_{r+1}\}}(k^s),$ let $a$ denote the number of $d_i$ equal to $m_i/2.$
\begEnu
\item $|\tStab_{J_T}(X)| = 2^a.$
\item $|L^{f,T}_{\{d_1,\ldots,d_{r+1}\}}(k^s)| = 2^r/2^a.$
\endEnu
\end{theorem}

Theorem \ref{thm:singular} is the second statement and Theorem \ref{thm:simpletransitivesingular} follows because the size of each orbit is $$|J_T|/|\tStab_{J_T}(X)|=2^r/2^a = |L^{f,T}_{\{d_1,\ldots,d_{r+1}\}}(k^s)|.$$ We will prove Theorem \ref{thm:regularcool} via a series of reductions.

\subsubsection*{Reduction on $d_1,\ldots,d_{r+1}$}

Suppose $X\in L^{f,T}_{\{d_1,\ldots,d_{r+1}\}}(K)$ with $d_i\geq1.$ Let $v_i$ denote an eigenvector of $T$ corresponding to $\al_i$. Since $T$ is regular, $v_i$ is unique up to scaling. The assumption $d_i\geq1$ then implies $v_i\in X.$ Let $H_i$ denote the hyperplane $v_i^\perp,$ and let $b$ denote the bilinear form associated to $Q$. Note $v_i\in H_i$ since there exists some $v'_i$ such that $(T-\al_i)v'_i=v_i,$ and hence
$$b(v_i,v_i)=b(v_i,(T-\al_i)v'_i)=b((T - \al_i)v_i,v'_i)=0.$$
For any $w\in H_i,$
$$b(v_i,Tw)=b(Tv_i,w) = b(\al_iv_i,w)=0.$$
Therefore, $T$ descends to a linear map $$\bar{T}_i:H_i/v_i\rightarrow H_i/v_i=:\bar{U}_i.$$ The quadratic form $Q$ descends to a non-degenerate quadratic form $\bar{Q}_i$ with respect to which $\bar{T}_i$ is regular self-adjoint with characteristic polynomial $f(x)/(x-\al_i)^2.$ Note this reduction can be described projectively as intersecting the quadric defined by $Q$ with the tangent plane to $v_i$, then projecting away from $v_i$.

Since $v_i\in X$ and $X$ is isotropic, we see $X\subset H_i.$ Let $\bar{X}_i$ denote the image of $X$ in $\bar{U}_i.$ It is immediate from the definition that $\bar{X}_i$ is $(n-1)$-dimensional, satisfying $$\bar{X}_i\subset\bar{X}_i^{\perp_{\bar{Q}_i}},\bar{T}_i\bar{X}_i\subset\bar{X}_i^{\perp_{\bar{Q}_i}},$$ and the maximal dimensions of  $\bar{T}_i$-stable subspaces in its intersection with the generalized eigenspaces are $d_1,\ldots,d_i-1,\ldots,d_{r+1}.$ We denote this reduction step by
$$\delta:L^{f,T}_{\{d_1,\ldots,d_{r+1}\}}(k^s) \xrightarrow{\,\,\sim\,\,} L^{f/(x-\al_i)^2,\bar{T}_i}_{\{d_1,\ldots,d_i-1,\ldots,d_{r+1}\}}(k^s).$$
$\delta$ is bijective, its inverse is given by taking the pre-image of the projection map $H_i\rightarrow \bar{U}_i.$

How are the stabilizers affected by this reduction? If $h(x)$ is any polynomial in $k^s[x],$ then $\delta(h(T)X)=h(\bar{T}_i)\bar{X}_i.$ Since $\delta$ is bijective, we conclude that $h(T)$ stabilizes $X$ if and only if $h(\bar{T}_i)$ stabilizes $\bar{X}_i.$ Note if $m_i\geq3,$ then $$h_i^T(x)=\frac{f(x)}{(x-\al_i)^{m_i}} = \frac{f(x)/(x-\al_i)^2}{(x-\al_i)^{m_i-2}}=h_i^{\bar{T}_i}(x).$$ Hence according to the explicit description given in Remark \ref{rmk:J2inPOregular}, $$h(T)\in J_T \Longleftrightarrow h(\bar{T}_i)\in J_{\bar{T}_i},\quad\mbox{hence}\quad|\tStab_{J_T}(X)| = |\tStab_{J_{\bar{T}_i}}(\bar{X}_i)|.$$

When $m_i=2,d_i=1,$ $\al_i$ is no longer an eigenvalue for $\bar{T}_i$. In this case, $$J_T = \lrg{h(T),1-2h_i(T)/h_i(\al_i)|h(\bar{T}_i)\in J_{\bar{T}_i}}.$$ Let $v'_i$ denote an element in $U_{i,T}$ such that $(T-\al_i)v'_i=v_i.$ Then $$U_{i,T}=\tSpan\{v_i,v'_i\},\quad\mbox{and}\quad b(v_i,v'_i)\neq0.$$ Since $v_i\in X$ and $X$ is isotropic, we see $$X = \tSpan\{v_i, X\cap\tSpan\{U_{j,T}\}_{j\neq i}\}.$$ Now $1-2h_i(T)/h_i(\al_i)$ sends $v_i$ to $-v_i$ and fixes every element in $\tSpan\{U_{j,T}\}_{j\neq i}.$ Therefore it stabilizes $X$ and hence $$|\tStab_{J_T}(X)| = 2|\tStab_{J_{\bar{T}_i}}(\bar{X}_i)|.$$ Note this case is precisely when $a$ decreases by 1 in this reduction step.

We summerize this reduction step in the following lemma.

\begin{lemma}\label{lem:d_reduction}
Suppose $d_i\geq1,$ then there is a bijection
$$\delta:L^{f,T}_{\{d_1,\ldots,d_{r+1}\}}(k^s) \xrightarrow{\,\,\sim\,\,} L^{f/(x-\al_i)^2,\bar{T}_i}_{\{d_1,\ldots,d_i-1,\ldots,d_{r+1}\}}(k^s).$$ The sizes of the stabilizers do not change, unless $m_i=2,d_i=1$ in which case it decreases by a factor of 2.
\end{lemma}

\subsubsection*{Reduction on $f$}

By the above reduction step, it remains to study $L^{f,T}_{\{0,0,\ldots,0\}}(k^s).$ We will describe the reduction map, state the corresponding result, then give the proof. However, since the proof is just hardcore linear algebra, we recommend the interested reader to prove it himself.

Suppose $\al$ is a root of $f$ of multiplicity $m\geq2.$ Let $X\in L^{f,T}_{\{0,0,\ldots,0\}}(k^s)$ be arbitrary. Let $v$ denote an eigenvector of $T$ with eigenvalue $\al.$ Suppose $v'\in U$ such that $(T-\al)v'=v.$ Since $$b(v,v)=b(v,(T-\al)v')=b((T-\al)v,v')=0,$$ we can consider the descent to $\bar{U}=v^\perp/v.$ As in the above reduction step, $Q$ descends to a non-degenerate quadratic form $\bar{Q}$ on $\bar{U}$ and $T$ descends to a regular self-adjoint operator $\bar{T}$ on $\bar{U}$ with characteristic polynomial $f(x)/(x-\al)^2.$

Observe that $v\notin X$ since $X$ contains no $T$-stable subspace. Therefore the map $U\rightarrow U/v$ is bijective when restricted to $X$. Consequently, $X\nsubseteq v^\perp$, for if otherwise the $(2n-1)$-dimensional vector space $v^\perp/v$ contains an $n$-dimensional isotropic subspace which is impossible. Now $X\cap v^\perp$ has dimension $n-1$ and we denote its bijective image in $v^\perp/v$ by $\bar{X}.$

\begin{lemma}\label{lem:f_reduction}
The above map sending $X$ to $\bar{X}$ defines a surjection
$$L^{f,T}_{\{0,0,\ldots,0\}}(k^s) \xrightarrow{\quad} L^{f/(x-\al)^2,\bar{T}}_{\{0,0,\ldots,0\}}(k^s).$$
This map is bijective if $m>2$ and is two-to-one if $m=2.$ In both cases, $$|\tStab_{J_T}(X)| = |\tStab_{J_{\bar{T}}}(\bar{X})|,\qquad \mbox{for any } X\in L^{f,T}_{\{0,\ldots,0\}}(K).$$
\end{lemma}

\bpr It is clear that $\bar{X}$ satisfies $\bar{X}\subset\bar{X}^\perp,\bar{T}\bar{X}\subset\bar{X}^\perp.$ If $\bar{X}$ contains a $\bar{T}$-stable subspace, then it must contain $v'+<v>.$ Hence $v' + cv\in X$ for some $c\in k^s.$ Then $v=(T-\al)(v' + cv)\in X^\perp$ contradicting $X\nsubseteq v^\perp.$ Therefore, $\bar{X}\in L^{f/(x-\al)^2,\bar{T}}_{\{0,0,\ldots,0\}}(k^s).$ We first prove surjectivity. Suppose $\bar{X}\in L^{f/(x-\al)^2,\bar{T}}_{\{0,0,\ldots,0\}}(k^s).$ Let $b_\al$ denote the bilinear form $$b_\al(u,u')=b(u,(T - \al)u').$$
Since $v$ lies in the kernel of $b_\al$, we see that $b_\al$ descends to a non-degenerate bilinear form on the $2n$ dimensional vector space $U/v.$ Denote by $\perp_\al$ the perpendicular space with respect to $b_\al.$ Since $\bar{X}$ is $n-1$ dimensional, $b_\al$ further descends to a non-degenerate bilinear form on the $2$-dimensional vector space $\bar{X}^{\perp_\al}/\bar{X}.$ It has two 1-dimensional isotropic lines, denote by $\bar{X}_1,\bar{X}_2$ their pre-images in $\bar{X}^{\perp_\al}.$

Suppose $m\geq3,$ let $v''$ be an element of $U$ such that $(T-\al)v''=v'.$ Then $$b_\al(v',v')=b(v',v)=b((T-\al)v'',v)=b(v'',(T-\al)v)=0.$$ Hence we might assume without loss of generality that $\bar{X}_1=\tSpan\{v'+<v>,\bar{X}\}\subset v^\perp/v.$ Since $\tSpan\{\bar{X}_1,\bar{X}_2\}$ has dimension $n+1,$ it is not isotropic with respect to $b_\al.$ Therefore, $b_\al(w,v')=b(w,v)\neq0$ for some $w+<v>\in\bar{X}_2.$ Up to scaling, we may assume $b(w,v)=1$ and by replacing $w$ by $w-\frac{1}{2}b(w,w)v,$ we may also assume $b(w,w)=0.$ Consider
\begin{eqnarray*}
X^w &=& \tSpan\{w,u - b(w,u)v\}_{u + <v>\in \bar{X}}\subset U,\\
(T-\al)X^w &=& \tSpan\{(T - \al)w, (T - \al)v\}.
\end{eqnarray*}
It is clear that $X^w\subset X^{w\perp}$ and $TX^w\subset X^{w\perp}$ with respect to $b$ by the construction of $w$. Since $w\notin v^\perp,$ we see $\bar{X^w}=\bar{X}.$ Since $b(w, c_2v)=c_2$, $X^w$ contains no non-zero vector of the form $c_2v$ and hence $X^w$ contains no non-zero $T$-stable subspace. We have now proved surjectivity when $m\geq3.$

Suppose now $X'\in L^{f,T}_{\{0,\ldots,0\}}(k^s)$ maps to $\bar{X}$. Then the image of $X'$ in $U/v$, denoted suggestively by $\bar{X'}_2$ is an $n$-plane isotropic to $b_\al,$ it contains $\bar{X}$ and is $b_\al$-orthogonal to $\bar{X}.$ Since it does not contain $v'+<v>,$ we conclude that $\bar{X'}_2=\bar{X}_2.$ Since the process from $\bar{X}_2$ to $X^w$ is just adjusting by adding the correct multiples of $v$, we see that $X'=X^w.$

Just as in the previous reduction step, when $m\geq3$, $J_T$ and $J_{\bar{T}}$ are represented by the same set of polynomials. It is clear that if $g(T)$ stabilizes $X$, then $g(\bar{T})$ stabilizes $\bar{X}.$ Conversely, if $g(\bar{T})$ stabilizes $\bar{X},$ then $g(T)$ sends $X$ to another $n$-plane that also maps to $\bar{X}.$ Since there is only one such $n$-plane, we conclude that $g(T)$ also stabilizes $X$. Therefore $$|\tStab_{J_T}(X)| = |\tStab_{J_{\bar{T}}}(\bar{X})|.$$

We now deal with the case $m=2$. Write $\bar{X}_1=\tSpan\{w_1+<v>,\bar{X}\}$ and $\bar{X}_2 = \tSpan\{w_2+<v>, \bar{X}\}.$ We claim $w_1\notin v^\perp$ and likewise same with $w_2.$ If for a contradiction that $w_1\in v^\perp,$ then $\bar{X}_1\subset v^\perp/v.$ When $m=2$, $v^\perp/v$ is the orthogonal (with respect to $b$) direct sum of all the generalized eigenspaces not containing $v,v'.$ Since $(T-\al)$ acts invertibly on generalized eigenspaces not containing $v,v'$, we see that $b_\al$ descends to a non-degenerate bilinear form on $v^\perp/v.$ However, $\bar{X}_1$ is isotropic of dimension $n$ and $v^\perp/v$ has dimension $2n-1$. Contradiction.

Finally, we lift each $\bar{X}_i$ to $X^{w_i}$ by adding an appropriate multiples of $v.$ The resulting $X^{w_i}$ both map to $\bar{X}$ under the reduction map. They are different from each other since their images in $U/v$ are different. Therefore we have proved surjectivity. The same argument as the above shows that $X^{w_1}$ and $X^{w_2}$ are precisely the two pre-images of $\bar{X}.$

Regarding stabilizers, we are in the situation where compared to $J_{\bar{T},}$ $J_T$ has an extra generator $h_0(T)=1-2h(T)/h(\al)$ where $h(x)=f(x)/(x-\al)^2.$ This extra generator fixes $v$ and acts as $-1$ on all the other generalized eigenspaces. Therefore $h_0(\bar{T})\bar{X}=\bar{X}$ and a simple computation shows that it switches $\bar{X}_1$ and $\bar{X_2}.$ If $g(\bar{T})$ stabilizes $\bar{X},$ then $g(T)$ either stabilizes $X^{w,1}$ or it sends $X^{w_1}$ to $X^{w_2}$, in which case $g(T)h_0(T)$ stabilizes $X^{w_1}.$ Therefore, the size of the stabilizers remains unchanged. \epr

\begin{cor}\label{cor:Lg0}
$|L^{f,T}_{\{0,0,\ldots,0\}}(k^s)| = 2^r$ and every element has trivial stabilizer in $J_T.$
\end{cor}

\bpr This follows from induction on the degree of $f$ and the classical result on generic intersection in odd dimension recalled in Section \ref{sec:povnonsin}. We write out the proof slightly differently from an induction argument so we can point out the differences between the contributions coming from roots of $f$ with odd multiplicity and the contributions from roots with even multiplicity. Rewrite the factorization of $f(x)$ as $$f(x) = \prod_{i=1}^{s_1+1}(x - \be_i)^{2n_i+1}\prod_{j=1}^{s_2}(x - \be'_j)^{2n'_j},$$ where each $\be_i$ is a root of $f(x)$ of odd multiplicity and each $\be'_j$ is a root of even multiplicity. Since $f(x)$ has odd degree, we know $s_1\geq 0$ and $s_1+s_2=r.$ Applying Lemma \ref{lem:f_reduction} repeatedly, one gets the following sequence of maps,
$$L^{f,T}_{\{0,0,\ldots,0\}}(k^s) \xrightarrow{\text{ 1 to 1 }} L^{\prod_i (x-\be_i)\cdot\prod_j (x-\be'_j)^2,T'}_{\{0,0,\ldots,0\}}(k^s) \xrightarrow{2^{s_2}\text{ to 1}} L^{\prod_i (x-\be_i),T''}_{\{0,0,\ldots,0\}}(k^s).$$
The last set has $2^{s_1}$ elements all of whose stabilizers are trivial. Applying Lemma \ref{lem:f_reduction} again, one concludes that every element in $L^{f,T}_{\{0,0,\ldots,0\}}(k^s)$ has trivial stabilizer as well. The above diagram shows that $|L^{f,T}_{\{0,0,\ldots,0\}}(k^s)|=2^{s_1+s_2}=2^r.$\epr

\textbf{Proof of Theorem \ref{thm:regularcool}: }Applying Lemma \ref{lem:d_reduction} repeatedly gives a bijection
$$L^{f,T}_{\{d_1,\ldots,d_{r+1}\}}(k^s)\xrightarrow[\delta]{\,\,\sim\,\,} L^{\prod_i(x-\al_i)^{m_i-2d_i},T'}_{\{0,0,\ldots,0\}}(k^s),$$ and for any $X\in L^{f,T}_{\{d_1,\ldots,d_{r+1}\}}(K),$ $$|\tStab_{J_T}(X)| = 2^a|\tStab_{J_{T'}}(\delta(X))|.$$

The polynomial $g(x)=\prod_i(x-\al_i)^{m_i-2d_i}$ has $r+1-a$ distinct roots, hence applying Corollary \ref{cor:Lg0} to $g$ completes the proof.\epr

\subsection{Even dimension}

Suppose now $U$ has dimension $N=2n+2$ for $n\geq1$. As above, suppose $Q=Q_1$ is non-degenerate and denote by $T$ the associated self-adjoint operator on $U$. As in Section \ref{sec:torofJ}, let $C$ be the (possibly singular) hyperelliptic curve parameterizing the rulings in the pencil. It is isomorphic over $k$, not canonically, to the hyperelliptic curve defined by $$y^2 = (-1)^{n+1}\det(Q)\det(xI - T)=\disc(Q)\det(xI - T).$$ To give a point on $C$ is the same as giving a quadric in the pencil along with a choice of ruling. Let $\w{C}$ denote its normalization. The geometric genus $p_g$ of $C$ is defined to be the genus of $\w{C}$. Let $C^{sm}$ denote the smooth locus of $C$.

\begin{lemma}\label{lem:maxdim1}
If $W$ is an $n+1$ dimensional subspace of $U\otimes k^s$ isotropic with respect to $Q_1,Q_2,$ then $W$ is $T$-stable, where $n\geq0$.
\end{lemma}

\bpr Take any $\lambda\in k$ that is not an eigenvalue of $T$. Then $W = W^{\perp_Q} = W^{\perp_{Q_\lambda}}.$
Hence, for any $w\in W,$ $(T-\lambda)w\in W^{\perp_{Q}}=W.$ In other words, $W$ is $T$-stable.\epr

\begin{proposition}\label{prop:maxdim2}
The base locus $B$ contains no $\bbp^n$ if and only if $p_g\geq0.$ When $C$ is reducible, or equivalently $p_g=-1,$ the base locus $B$ contains a unique $\bbp^n.$
\end{proposition}

\bpr Without loss of generality, assume $k$ is separably closed. Suppose $W$ is an $n+1$ dimensional subspace of $U$ such that $\bbp W\subset B.$ The above lemma says $W$ contains an eigenvector $v$ of $T$. Since $W$ is isotropic, the eigenvalue of $v$ has multiplicity at least 2. One can now reduce the problem to $\bar{U}=v^\perp/v$ and $\bar{W}$ is $n$-dimensional. Applying the above lemma and reduction repeatedly until $\dim\bar{U}=2$ and $\dim\bar{W}=1.$ Apply the above lemma again, we see that $\bar{T}$ has a repeated eigenvalue and hence all the generalized eigenspaces of $T$ have even dimension which implies that $C$ is reducible. Conversely when $C$ is reducible, $\bar{W}$ is the unique 1-dimensional eigenspace of $\bar{T}$ hence proving uniqueness. Existence follows from running the argument backwards.\epr

Let $F_0$ denote the following variety over $k$,
$$F_0 = \{\pX|\dim \pX = n-1, \pX\subset B\}.$$
In view of the above subsection and Example \ref{ex:nonWeier}, we impose an open condition and look at the following variety,
\begin{equation}\label{eq:defF}
F = \{\pX\in F_0|\tSpan\{X,TX\}\mbox{ has no non-zero }T\mbox{-stable subspace}\}.
\end{equation}

\begin{lemma}\label{lem:equivF}
Suppose $p_g\geq0,$ then
\begin{eqnarray*}
F &=& \{\pX\in F_0|X\nsubseteq v^\perp,\mbox{ for all singular points }[v]\in B\}\\
&=& \{\pX\in F_0|[v]\notin\pX,\mbox{ for all singular points }[v]\in B\}.
\end{eqnarray*}
\end{lemma}

\bpr Suppose $\pX\in F$. Let $[v]$ be any singular point of $B$, since $v$ is an eigenvector, $v\notin X.$ If $X\subset v^\perp,$ then $\bbp(\tSpan\{X,v\})$ is a $\bbp^n$ contained in $B$, contradicting Proposition \ref{prop:maxdim2}.

Conversely, suppose $\pX\notin F,$ then $v\in\tSpan\{X,TX\}$ for some eigenvector $v$ of $T$. Since $X$ is a isotropic with respect to every quadric in the pencil, we see that $v\in\tSpan\{X,TX\}\subset X^\perp$ and hence $X\subset v^\perp.$

For the second equality, suppose first $X\subset v^\perp$ for some singular $[v]\in B.$ If $v\notin X$, then aftering reduction to $v^\perp/v,$ $(X\cap v^\perp)/v$ has dimension $n$ which contradicts Proposition \ref{prop:maxdim2}. Hence $v\in X.$ Conversely, if $v\in X,$ then $X\subset v^\perp$ as above.\epr

\begin{remark} The main reason why $F$ was defined as in \eqref{eq:defF} instead of the more conceptual ones in Lemma \ref{lem:equivF} is that there is still some interesting geometry when $p_g=-1$ as we will see towards the end of the paper, and in that case, \eqref{eq:defF} is the more appropriate definition.
\end{remark}

\begin{theorem}\label{thm:MainRegEven2}
Suppose $p_g\geq0$ and $C$ only has nodal singularities. Then there is a commutative algebraic group structure $+_G$ defined over $k$ on the disconnected variety $$G = \picz(C) \dcup F \dcup \pico(C) \dcup F'$$ such that,
\begEnu
\item $G^0 = \picz(C)$ with component group $G/G^0 \simeq \bbz/4,$
\item $F'$ is isomorphic to $F$ as varieties via the inversion map $-1_G$,
\item the group law extends that on $H = \Pic{C}/D_0 \simeq \picz(C) \dcup \pico(C)$ where $D_0$ is the hyperelliptic class.
\endEnu
\end{theorem}

From now on, we assume that $p_g\geq0.$ Since the base locus contains no $\bbp^n$, one can define $\tau:C\times F_0\rightarrow F_0$ as in the generic case.

\begin{lemma}\label{lem:taugen} $\tau$ restricts to a morphism $C^{sm}\times F\rightarrow F.$
\end{lemma}

\bpr Recall that given a pair $(c,\pX)\in C^{sm}\times F,$ there is a unique $\pY\simeq\bbp^n$ in the quadric and the ruling defined by $c$, then $\tau(c,\pX)$ is the residual intersection of $\pY$ with the base locus. The claim here is that $\tau(c,\pX)\in F.$ Suppose for a contradiction that $\pX':=\tau(c,\pX)\in F_0-F.$ Then by Lemma \ref{lem:equivF}, there exists a singular point $[v]\in B$ such that $X'\subset v^\perp.$ Hence the linear space $\tSpan\{X',v\}$ is isotropic with respect to every quadric in the pencil. Proposition \ref{prop:maxdim2} implies that $v\in X'.$ Since $X$ and $X'$ intersect at codimension 1 and $v\notin X,$ we see that $$\pY=\tSpan\{\pX,\tau(c,\pX)\}=\tSpan\{\pX,[v]\}.$$ Since $\pY$ lies in the quadric $Q_\al$ where $\al$ is the eigenvalue of $v,$ we see that $c=(\al,0)\notin C^{sm}.$ Contradiction.\epr

As in the generic case, one obtains an action of $C^{sm}$ on $F\dcup F',$
\begin{equation}\label{eq:actionofCgen}
\pX + (c) = -\tau(\bar{c})\pX,\quad\quad -\pX + (c) = \tau(c)\pX.
\end{equation}
This action extends to an action of $\tDiv(C^{sm})$ on $F\dcup F'.$ To show that this descends to a simply-transitive action of $\picz(C),$ we assume $k=k^a$ and work over the algebraic closure. Let $v$ be an eigenvector with eigenvalue $\al$ of multiplicity $m\geq2$. As usual, let $(\bar{U},\bar{Q})$ denote the $2n$-dimensional quadratic space $v^\perp/v,$ let $\bar{T}$ denote the descent of $T$ to $\bar{U}$. Let $\bar{C}$ denote the (possible singular) hyperelliptic curve
$$y^2 = \disc(\bar{Q})\det(xI - \bar{T}) = \disc(Q)\det(xI - T)/(x-\al)^2.$$
Note $\bar{C}\rightarrow C$ is a partial normalization of $C$. There is a natural inclusion $\iota:C^{sm}\hookrightarrow\bar{C}^{sm}.$ Define $\bar{F}$ and $\bar{F}_0$ in the analogous way as $F$ and $F_0$. Suppose $\pX\in F,$ write $\bar{X}=(X\cap v^\perp)/v.$ Lemma \ref{lem:equivF} implies that $\bar{X}$ has the correct dimension. It is clear therefore $\bar{X}\in\bar{F}_0$.

\begin{lemma}\label{lem:reductionwell}
$\tSpan\{\bar{X},\bar{T}\bar{X}\}$ has no non-zero $\bar{T}$-stable subspace.
\end{lemma}

\bpr Note this is immediate when $C$ has only nodal singularities for this reduction step kills the $\al$-generalized eigenspace and leaves the rest unchanged. In general, by Lemma \ref{lem:equivF}, it suffices to show $\bar{X}$ does not contain any singular point of $\bar{B}$. Let $v'\in U$ be such that $(T-\al)v'=v.$ Then $\bar{X}$ could possibly contain a singular point of $\bar{B}$ if $m\geq4$ and $v'+cv\in X$ for some $c\in k.$ The latter condition implies $v=(T-\al)(v'+cv)\in X^\perp$ contradicting $X\nsubseteq v^\perp.$\epr

Denote this reduction step by $\delta_v:F\rightarrow \bar{F}.$ We now have the following commutative diagram,
\begin{displaymath}
\xymatrix{
C^{sm}\times F\ar[d]^{\iota\times\delta_v}\ar[r]&F\ar[d]^{\delta_v}\\
\bar{C}^{sm}\times\bar{F}\ar[r]&\bar{F}}
\end{displaymath}
The natural map $\bar{C}\rightarrow C$ induces a map $J(C)\rightarrow J(\bar{C})$ on their Jacobians with kernel either $\bbg_m$ if the multiplicity $m$ of $\al$ is 2, or $\bbg_a$ if $m\geq3.$ We now show that $\delta_v$ is surjective and the preimage of every point is isomorphic to $\ker(J(C)\rightarrow J(\bar{C})).$ Let $b_\al$ denote the bilinear form $b_\al(u,u')=b(u,(T-\al)u')$ and by $\perp_\al$ the operation of taking perpendicular space with respect to $b_\al$. Fix any $\bar{X}\in\bar{F}$. The bilinear form $b_\al$ descends to a non-degenerate form on the $2n+1$ dimensional space $U/v.$ Inside this space, we have
\begin{eqnarray*}
\dim\, \bar{X}^{\perp_\al}/\bar{X} &=& 3,\\
\dim\, (\bar{X}^{\perp_\al} \cap v^\perp)/\bar{X} &=&2.
\end{eqnarray*}
Stated in a different way, $b_\al$ defines a smooth conic $C_0$ in $\bbp^2 = \bbp(\bar{X}^{\perp_\al}/\bar{X})$ and $l=\bbp((\bar{X}^{\perp_\al} \cap v^\perp)/\bar{X})$ is a line intersecting the conic at either one point or two points.
\begin{figure}[h]
                \centering
                \includegraphics[width=4in]{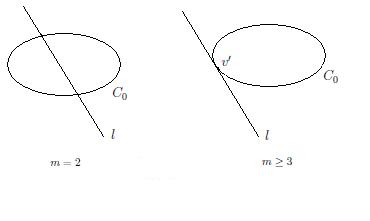}
                \label{fig:mequals2}
\end{figure}

\begin{lemma} $l$ intersects $C_0$ tangentially if and only if $m\geq3,$ in which case the point of intersection is $[v'+<v>+\bar{X}],$ where $v'\in U$ is such that $(T - \al)v'=v.$
\end{lemma}

\bpr Suppose $l$ intersects $C_0$ at a point $w+<v>+\bar{X}$. To say $l$ intersects $C_0$ tangentially at $w+<v>+\bar{X}$ is equivalent to saying
\begin{equation}\label{eq:condtang}
w+<v>\in\bar{X}^{\perp_\al},\quad b_\al(w,w)=0, \quad w^{\perp_\al}\cap\bar{X}^{\perp_\al}=v^{\perp}\cap\bar{X}^{\perp_\al}.
\end{equation}
Since $v'^{\perp_\al}=v^\perp,$ we have $v'\in (v^\perp)^{\perp_\al}.$ Thus $(v^{\perp}\cap\bar{X}^{\perp_\al})^{\perp_\al}\cap \bar{X}^{\perp_\al}$ is the line spanned by $v'+<v>.$ Since $w\in (w^{\perp_\al})^{\perp_\al}$, we see that up to scaling $w+<v> = v'+<v>$. Finally, $b_\al(v',v')=b(v',v)=0$ if and only if $m\geq3.$

Conversely, suppose $m\geq3,$ then $v'\in v^\perp$ and it is easy to see $w=v'$ satisfies \eqref{eq:condtang}.\epr

Now given any point $[w+<v>+\bar{X}]\in C_0 - l,$ we can proceed to find a lift of $\bar{X}$ to $X^w\in F$ as follows. Since $b(w,v)\neq0,$ we can choose a lift of $w\in U$ unique up to scaling such that $b(w,w)=0$ by adding an appropriate multiple of $v$, then take
$$X^w = \tSpan\{w,u - \frac{b(w,u)}{b(w,v)}v\}_{u + <v>\in \bar{X}}\subset U.$$
To check $X^w\in F,$ we only need to check $X^w\nsubseteq v^\perp$, which is clear since $b(w,v)\neq0.$ For any two points in $C_0-l,$ the corresponding lifts to $F$ are distinct as they have different images in $U/v.$ Lastly, if $X\in F$ such that $\delta_v(X)=\bar{X},$ then the image of $X$ in $U/v$ must be of the form $\tSpan\{\bar{X},w+<v>\}$ for some $w+<v>+\bar{X}\in C_0-l.$ Therefore, we have prove the following proposition.

\begin{proposition}\label{prop:onestepreduction}
$\delta_v:F\to\bar{F}$ is surjective. The fibers are isomorphic to either (a conic minus a point)$\simeq\bbg_a$ when $m\geq3$, or (a conic minus two points)$\simeq\bbg_m$ when $m=2$. The kernel of the map $J(C)\to J(\bar{C})$ has the same property.
\end{proposition}

One can now apply this reduction with any singular point of $\bar{B}$ and so on. For each $i$ such that $m_i\geq 2,$ let $v_{i,1}$ denote an eigenvector of $T$ with eigenvalue $m_i,$ and let $v_{i,j}$ be such that $(T-\al_i)v_{i,j}=v_{i,j-1}$ for $j=2,\ldots,\lfloor\frac{m_i-1}{2}\rfloor.$ Let $V$ denote the linear span of all such $v_{i,j}.$ The above reduction will terminate at the $2p_g+2$ dimensional vector space $\w{U}=V^\perp/V$. The data $(Q,T)$ descends to $(\w{Q},\w{T})$ on the $2p_g+2$ dimensional vector space $\w{U}=V^\perp/V$ with $\w{T}$ regular semi-simple. Let $\w{F}$ denote the variety of (linear) $p_g$-dimensional common isotropic subspaces $\w{X}\subset\w{U}.$ Let $\delta:F\rightarrow \w{F}$ denote the composite of all the reductions. The associated smooth hyperellipitic curve $\w{C}$ is the normalization of $C$. Note that if $k$ is arbitrary, then $V$ is defined over $k$ and the composite $\delta$ is defined over $k$. We summarize the above discussion into the following Theorem.

\begin{theorem}\label{thm:totalreduction} Suppose $p_g\geq0$ and $k$ is algebraically closed. Then:
\begEnu
\item The map $\delta:F\onto \w{F}$ is surjective. The pre-image of every point has a filtration with $\bbg_a$ and $\bbg_m$ factors. The kernel of the natural map $J(C)\to J(\w{C})$ has a filtration with the same factors.
\item There is an action of $\mbox{Div}^0(C^{sm})$ on $F$ that descends to the simply-transitive action of $J(\w{C})$ on $\w{F}.$
\endEnu
\end{theorem}

Therefore to prove Theorem \ref{thm:MainRegEven2}, it remains to show that the action of $\tDiv(C^{sm})$ on $F\dcup F'$ descends to a simply-transitive action of $\text{Pic}(C)$ on $F\dcup F'$. Once again we pass to the algebraic closure and use the same formal argument as in the generic case. We list the ``non-formal'' results one needs to verify in the regular case.

\begEnu
\item Lemma \ref{lem:Don32}, which allows one to define the $\infty$-minimal form of a divisor class $[D]\in J(C)$ and hence a morphism $\vphi:J\rightarrow \text{Aut}(F).$ Here we need to assume that $C$ has a smooth Weierstrass point.
\item Show $\vphi$ is a group homomorphism, to conclude that principal divisors supported on $C^{sm}$ act trivially on $F\dcup F'$.
\item The existence part of Lemma \ref{lem:Don26}, to conclude that the action of $J$ on $F$ is transitive.
\item The uniqueness part of Lemma \ref{lem:Don26}, to conclude that the action is simply-transitive.
\endEnu

Lemma \ref{lem:Don32} still holds in the singular case because Riemann-Roch holds in the singular case (\cite{Hartshorne}). Suppose $C$ has a smooth Weierstrass point $\infty$, which it always has if $C$ only has nodal singularity and $p_g\geq0.$ Every class $[D]\in J(C)$ has a $\infty$-minimal form $[D'-r(\infty)]$ where $D'$ is effective of degree $r\geq n$ supported on $C^{sm}$ and $h^0(D)=1$. This allows us to define a morphism of varieties $\varphi:J\rightarrow \text{Aut}(F).$ The image of $\vphi$ lies in a commutative subvariety of $\text{Aut}(F)$.

We now specialize to the case where $C$ only has nodal singularities, so $J$ is an extension of an abelian variety $\w{J}$ of dimension $p_g$ by an $n-p_g$ dimensional torus $S$.

\begin{lemma}\label{lem:rigidity}
$\varphi$ is a morphism of algebraic groups.
\end{lemma}

\bpr The proof is very similar to the proof that a morphism between semi-abelian varieties mapping the identity to the identity is a group homomorphism. For any $s\in S,$ its image in $\w{J}$ is 0, hence it acts on the fibers of the map $\delta:F\rightarrow\w{F}$ which are also tori. Therefore $\varphi_{|S}$ is a group homomorphism. For any $a\in J,$ we define
$\varphi_a:S\rightarrow \text{Aut}(F)$ by $$\varphi_a(s) = \vphi(a)\vphi(s)\vphi(as)^{-1}.$$ Fix any $x\in F$, we have $\delta(x) = \delta(\vphi_a(s)(x))$. Let $S'$ denote the fiber of $\delta$ over $\delta(x),$ we have thus defined a map $\vphi_{a,x}:S\rightarrow S'$ between tori, which is automatically a group homomorphism. Letting $a$ vary, one obtains a map $\vphi_x:J\rightarrow\text{End}(S,S').$ Since $J$ is connected and $\text{End}(S,S')$ is discrete, $\vphi_x$ is constant. Taking any $s\in S,$ we see $\vphi_x(a)=\vphi_x(s)$ is the trivial map $S\rightarrow S'.$ Letting $x$ vary, we have proved that
\begin{equation}\label{eq:semiab}
\vphi(a)\vphi(s)=\vphi(as),\quad\forall a\in J,s\in S.
\end{equation}
Now fix $a\in J$ and view $\vphi_a$ as a morphism $J\rightarrow \text{Aut}(F).$ Since $\vphi_a$ vanishes on $S$ and \eqref{eq:semiab} allows us to descend $\vphi_a$ to a morphism $\w{J}\rightarrow \text{Aut}(F)$. Once again, fixing any $x\in F,$ $\vphi_a(a')$ acts on the fiber over $\delta(x).$ Hence we have a morphism $\vphi_{a,x}:\w{J}\rightarrow S'$ which is trivial since $\w{J}$ is an abelian variety and $S'$ is a torus. Letting $x$ vary, one sees that $\vphi_a$ is trivial. Letting $a$ vary gives the desired result.\epr

As in the proof of Proposition \ref{prop:principaldies}, we have shown that principal divisors supported on $C^{sm}$ act trivially on $F\dcup F'$. Next we show transitivity of this action. Since $\tDiv(C^{sm})$ also acts on $F_0\dcup F'_0$ and $F\dcup F'$ is open in $F_0\dcup F'_0$, by taking Zariski closure one sees that principal divisors supported on $C^{sm}$ act trivially on $F_0\dcup F'_0$. Since being supported on $C^{sm}$ is also an open condition, one also has that principal divisors on $C$ act trivially on $F_0$. The existence part of Lemma \ref{lem:Don26} can be applied to $F_0$ since the defining map $C\rightarrow \bbp^1$ admits no section. In other words, given $x,x'\in F,$ view them as in $F_0$ where there exists an effective divisor $D\in \tDiv(C)$ such that $x + D = \pm x'.$ Let $D'$ be a divisor supported on $C^{sm}$ linearly equivalent to $D$. Since principal divisors on $C$ act trivially, $x+D'=x+D=\pm x'.$ Transitivity then follows from the formal argument in the proof of Proposition \ref{prop:transitive}. Note here the existence of a smooth Weierstrass point is needed because we need to know there exists $\pX\in F$ such that $T_{\pX}B\simeq \bbp^n$.

The uniqueness part of Lemma \ref{lem:Don26} also holds for $F_0$. The argument in \cite{Donagi} works since there is no injective map from $\bbp^1$ to $C$ when the arithmetic genus $n$ of $C$ is at least 1. The same formal argument in the generic case then implies that only principal divisors act trivially. Once again, the existence of a smooth Weierstrass point $\infty$ is also needed, for we need the analogous Example \ref{ex:Weiergen} to know there are finitely many element of $F$ fixed by $\tau(\infty).$ We have now finished the proof of Theorem \ref{thm:MainRegEven2}. The following result is immediate from Theorem \ref{thm:MainRegEven2} and Theorem \ref{thm:totalreduction}.

\begin{cor}\label{cor:SESofG}
Suppose $p_g\geq0$ and $C$ only has nodal singularities. Then the short exact sequence
$$1\rightarrow T\rightarrow J(C)\rightarrow J(\w{C})\rightarrow 1$$
extends to a short exact sequence
$$1\rightarrow T\rightarrow G\rightarrow \w{G}\rightarrow 1,$$
where $G = \picz(C) \dcup F \dcup \pico(C) \dcup F'$ and $\w{G} = \picz(\w{C}) \dcup \w{F} \dcup \pico(\w{C}) \dcup \w{F}'$ are the corresponding disconnected groups of four components.
\end{cor}

Now over the algebraic closure, after identifying $F$ with $J(C)$, one can obtain a compactification of $J(C)$ by taking $F_0.$ Recall for any singular $[v]\in B,$ we have the reduction map $\delta_v:F_0\rightarrow \bar{F}_0.$ Note this map might not be a morphism. The composition of all the reduction map gives $\delta:F_0\rightarrow \w{F}\simeq J(\w{C}).$ Each fiber of $\delta_v$ intersects $F_0\backslash F$ at one point, obtained by taking the preimage of $\pX\in \bar{F}_0$ under the map $v^\perp\rightarrow v^\perp/v.$

\begin{cor}\label{cor:compactJ} Suppose $p_g\geq0$ and $C$ has only nodal singularities, then $F_0$ is a compactification of $J(C)$ by adding one point to each $\bbg_m$ factor of the fiber over $J(\w{C})$. 
\end{cor}

We expect that the condition on $C$ having only nodal singularities is unnecessary. If Theorem \ref{thm:MainRegEven2} is proved without this condition, then Corollary \ref{cor:compactJ} also holds without this condition. The compactification $F_0$ is not smooth.

Since multiplication by 2 and 4 are still surjective on $J$, we can lift $F$ to a torsor of $J[4]$ by taking
$$F[4] = \{\pX\in F|\pX +_G \pX +_G \pX +_G \pX = 0\}.$$
When $\pico(C)(k)\neq\emptyset,$ $[F]$ is 2-torsion. For each $[D_1]\in\pico(C)(k),$ we obtain a lift of $F$ to a torsor of $J[2]$ by taking
$$F[2]_{[D_1]}=\{\pX\in F|\pX +_G \pX = [D_1]\}.$$

\subsubsection{Example: smooth rational Weierstrass point}
\begin{example}\label{ex:Weiergen}
Suppose $C$ has a smooth rational Weierstrass point. By moving this point to $\infty,$ we assume that $Q_1$ is singular with cone point $[v_\infty].$ Let $H=v_\infty^{\perp_{Q_2}}$ be the hyperplane in $U$ orthogonal to $v_\infty$ with respect to $Q_2$. Then $\tau(\infty)$ is induced by the linear map on $U$ that fixes $H$ and sends $v_\infty$ to $-v_\infty.$ Hence,
\begin{equation}\label{eq:F2weiergen}
F[2]_\infty = \{\pX\in F|\pX\subset B\cap \bbp H\}.
\end{equation}

Just as in the generic case, when restricted to the $2n+1$ dimensional vector space $H$, $Q_1$ and $Q_2$ span a regular pencil $\scl_H$. Moreover, $Q_{1|H}$ is non-degenerate and $T_{|H}$ restricts to the self-adjoint operator on $H$ associated to the pencil $\scl_H$ as defined in \eqref{eq:defT}. The right hand side of \eqref{eq:F2weiergen} is precisely $L^{f,T_{|H}}_{\{0,0,\ldots,0\}}$ as defined in the odd dimension case. Now $J[2]$ acts on $F[2]_\infty$ via the action of $J$ and on $L^{f,T_{|H}}_{\{0,0,\ldots,0\}}$ via the identification $J[2]\simeq \tStab_{\PO(H,Q_{1|H})}(T).$ As in the generic case, these two actions coincide.
\end{example}

\subsubsection{Example: rational non-Weierstrass point}
\begin{example}\label{ex:nonWeiergen}
Suppose $C$ has a rational non-Weierstrass point, or equivalently, $\scl$ has a rational quadric with discriminant 1. By moving the point to infinity, we assume that $Q_1$ has discriminant 1. Its two rulings are therefore defined over $k$. Let $Y_0$ denote one of the rulings and let $\infty\in C(k)$ denote the point corresponding to the quadric $Q_1$ and the ruling $Y_0$. Denote by $\infty'$ the conjugate of $\infty$ under the hyperelliptic involution. As in the generic case, we have,
$$F[2]_\infty = \{\pX\in F|\pX = \tau(\infty)\pX\} = \{\pX\in F | \tSpan\{\pX,\bbp(TX)\}\sim Y_0\}.$$
The latter condition means $\tSpan\{\pX,\bbp(TX)\}\simeq\bbp^n$ is contained in $Q_1$ in the ruling $Y_0$.

Fix now, the monic polynomial $f$ of degree $2n+2$ splitting completely over $k^s$, the quadratic form $Q=Q_1$ of discriminant 1, and for every field $k'$ containing $k$, define
$$V_f(k') = \{T:U\otimes k'\rightarrow U\otimes k'|T\mbox{ is self-adjoint and regular with minimal polynomial }f(x)\}.$$
For every field $k'$ containing $k$, and every $T\in V_f(k'),$ let $W_T(k')$ denote the set of (linear) $n$-dimensional $k'$-subspaces $X$ of $U\otimes k'$ such that $\tSpan\{X,TX\}\sim Y_0$. That is to say the linear space $\tSpan\{X,TX\}$ is an $(n+1)$-dimensional isotropic subspace with respect to $Q$ that lies inside the ruling $Y_0$ over $k'.$ As before, we define
$$W_f(k')=\{(T,X)|T\in V_f(k'),X\in W_T(k')\}.$$ There is a Galois invariant action of $\PSO(U,Q)=\SO(U,Q)/(\pm1)$ on $W_f:$ $$g.(T,X)=(gTg^{-1},gX).$$

Recall that $U\otimes k^s$ breaks up as the orthogonal direct sum of generalized eigenspaces $U_{i,T}$ of $T$ of dimension $m_i$. For any linear space $X$, we defined $\dim_{i,T}(X)$ to be the dimension of the maximal $T$-stable subspace in $(X\otimes k^s)\cap U_{i,T}$. For any sequence of integers $d_1,\ldots,d_{r+1}$ such that $0\leq d_i\leq m_i/2,$ we defined for any field $k'$ containing $k$,
$$L^{f,T}_{\{d_1,\ldots,d_{r+1}\}}(k')=\{X\simeq (k')^n|X\subset X^\perp,TX\subset X^\perp, \dim_{i,T}(X)= d_i\}.$$
In view of the definition of $F$, we define,
$$L'^{f,T}_{\{d_1,\ldots,d_{r+1}\}}(k')=\{X\in W_T(k')|\dim_{i,T}(\tSpan\{X,TX\})=d_i\},$$
$$W'^f_{\{d_1,\ldots,d_{r+1}\}}(k') = \{(T,X)|T\in V_f(k'),X\in L'^{f,T}_{\{d_1,\ldots,d_{r+1}\}}(k')\}.$$
We make no assumption on the reducibility of the associated hyperelliptic curve $C$ but assume instead a weaker condition,
\begin{equation}\label{eq:XisnotTstable}
d_1+\cdots+d_{r+1} < n + 1 = \dim\,\tSpan\{X,TX\}.
\end{equation}
This condition is equivalent to saying $\tSpan\{X,TX\}$ is not $T$-stable. Let $s_1$ denote the number of roots of $f$ with odd multiplicity. Then the maximum $d_1+\cdots+d_{r+1}$ could reach is $n+1-s_1/2.$ If \eqref{eq:XisnotTstable} fails, then we must have $s_1=0$ and hence $C$ is reducible. If one uses $L^{f,T}$ instead of $L'^{f,T}$ or if one does not assume \eqref{eq:XisnotTstable}, then there will be infinitely many choices for $X$ when $C$ is reducible. See Example \ref{ex:22} and Example \ref{ex:400}.

As one would expect from the odd case, the main theorem we are heading towards is the following:

\begin{theorem}\label{thm:singulareven}
Suppose $d_1+\cdots+d_{r+1} < n + 1$, then $|L'^{f,T}_{\{d_1,\ldots,d_{r+1}\}}(k^s)|=2^r/2^a,$ where $a$ is the number of $d_i$'s equal to $m_i/2.$
\end{theorem}

The action of $\PSO(U,Q)$ preserves the decomposition of $U\otimes k^s$ into generalized eigenspaces. Therefore one obtains a Galois invariant action of $\PSO(U,Q)$ on $W'^f_{\{d_1,\ldots,d_{r+1}\}}.$

\begin{theorem}\label{thm:simpletransitivesingulareven}
Suppose $d_1+\cdots+d_{r+1} < n + 1$, then $\PSO(U,Q)(k^s)$ acts on $W'^f_{\{d_1,\ldots,d_{r+1}\}}(k^s)$ simply-transitively if $a=0$ and transitively if $a>0.$
\end{theorem}

\begin{cor}\label{thm:simpletransitivesmoothregulareven} Suppose $d_1+\cdots+d_{r+1} < n + 1$. Then $\PSO(V,Q)(k')$ acts simply-transitively on $W'^f_{\{0,\ldots,0\}}(k')$ for any field $k'$ over $k$.
\end{cor}

\bpr Same descent argument as in the proof of Corollary \ref{thm:simpletransitivesmooth}.\epr

We begin by studying the conjugation action of $\PSO(U,Q)$ on $V_f.$

\begin{proposition}\label{prop:POtransitiveregulareven}
The action of $\PSO(U,Q)$ on $V_f$ has a unique geometric orbit. For any $T\in V_f(k')$ defined over some field $k'$ over $k$, its stabilizer scheme $\tStab(T)$ is isomorphic to $(\mbox{Res}_{L'/k'}\mu_2)_{N=1}/\mu_2\simeq J_{k'}[2]$ as group schemes over $k'$ where $L' = k'[x]/f(x).$ In particular, $\tStab_{PSO(U,Q)}(T)(k^s)$ is an elementary abelian 2-group of order $2^r.$
\end{proposition}

\bpr cf. Proposition \ref{prop:POtransitiveregular}.\epr

\begin{remark}\label{rmk:J2inPOregulareven}
A more explicit description for the stabilizer as polynomials in $T$ is almost identical to the odd case as given in Remark \ref{rmk:J2inPOregular}. For each $i=1,\ldots,r+1,$ define $h_i^T(x) = f(x)/(x-\al_i)^{m_i}.$ Then
\begin{eqnarray*}
\mu_2(K[T]^\times) &=& \left\{\prod_{i\in I}\left(1-2\frac{h_i^T(T)}{h_i^T(\al_i)}\right)
\right\}_{I\subset\{1,\ldots,r+1\},2\,|\,|I|}\\
&=&\left\{1-2\sum_{i\in I}\frac{h_i^T(T)}{h_i^T(\al_i)}\right\}_{I\subset\{1,\ldots,r+1\},2\,|\,|I|}.
\end{eqnarray*}
For any $I\subset\{1,\ldots,r+1\}$ and any $j\notin I,$ since $(x-\al_j)^{m_j}$ divides $h_i(x)$ in $K[x]$ and $(T-\al_j)^{m_j}$ kills all the generalized eigenspaces $U_{j,T}$, $$1-2\sum_{i\in I}\frac{h_i^T(T)}{h_i^T(\al_i)}$$ acts trivially on $U_{j,T}$.
\end{remark}

\begin{cor}
For any $T,T'\in V_f(K)$, there exists a bijection $$L'^{f,T}_{\{d_1,\ldots,d_{r+1}\}}(K)\longleftrightarrow L'^{f,T'}_{\{d_1,\ldots,d_{r+1}\}}(K).$$ 
\end{cor}

\bpr Suppose $g\in \PSO(U,Q)(k^s)$ conjugates $T$ to $T'$, then the left action by $g$ on $\mbox{Gr}(n,U)$ gives the desired bijection.\epr

Also by Proposition \ref{prop:POtransitiveregular}, for any $T\in V_f(K),$ its stabilizer $J_T$ in $\PSO(U,Q)(K)$ acts on $L'^{f,T}_{\{d_1,\ldots,d_{r+1}\}}(K)$. We rephrase the main theorems as follows.

\begin{theorem}\label{thm:regularcooleven}
Suppose $d_1+\cdots+d_{r+1} < n + 1$. For any $X \in L'^{f,T}_{\{d_1,\ldots,d_{r+1}\}}(k^s),$ let $a$ denote the number of $d_i$ equal to $m_i/2.$
\begEnu
\item $|\tStab_{J_T}(X)| = 2^a.$
\item $|L'^{f,T}_{\{d_1,\ldots,d_{r+1}\}}(k^s)| = 2^r/2^a.$
\endEnu
\end{theorem}

Theorem \ref{thm:singulareven} is the second statement and Theorem \ref{thm:simpletransitivesingulareven} follows because the size of each orbit is $$|J_T|/|\tStab_{J_T}(X)|=2^r/2^a = |L'^{f,T}_{\{d_1,\ldots,d_{r+1}\}}(k^s)|.$$ We will prove Theorem \ref{thm:regularcool} via a series of reductions.

One major difference from the odd case is that one should forget about the rulings in the following reductions. Namely, consider instead
$$W^*_T(k^s)=\{X\in\mbox{Gr}(n,U\otimes k^s)|\tSpan\{X,TX\}\mbox{ is }n+1\mbox{ dimensional and isotropic}\}.$$
Observe that $W^*_T(k^s)$ can be divided into two components, one of which is $W_T(k^s),$ corresponding to which ruling $\tSpan\{X,TX\}$ lies in. The two components are in bijection with each other via an element in $\tStab_{PO}(T)$ not in $\tStab_{PSO}(T).$ One defines similarly $L'^{f,T,*}_{\{d_1,\ldots,d_{r+1}\}}(k^s).$

\subsubsection*{Reduction on $d_1,\ldots,d_{r+1}$}

Suppose $X\in L'^{f,T,*}_{\{d_1,\ldots,d_{r+1}\}}(k^s)$ with $d_i\geq1.$ Let $v_i$ denote an eigenvector of $T$ corresponding to $\al_i$. Since $T$ is regular, $v_i$ is unique up to scaling. The assumption $d_i\geq1$ implies $v_i\in \tSpan\{X,TX\}.$ Hence $$X\subset\tSpan\{X,TX\}^\perp\subset v_i^\perp=:H_i.$$ Let $b$ denote the bilinear form associated to $Q$. As before, the data $(Q,T)$ descentds to $(\bar{Q}_i,\bar{T}_i)$ on $\bar{U}_i:=H_i/v_i$ and $\bar{T}_i$ is regular with characteristic polynomial $f(x)/(x-\al_i)^2.$ Let $\bar{X}_i$ denote the image of $X$ in $\bar{U}_i.$ Then $\tSpan\{\bar{X}_i,\bar{T}_i\bar{X}_i\}$ is an isotropic $n$-plane with respect to $\bar{Q}_i$, and the maximal dimensions of $\bar{T}_i$-stable subspaces in the intersection of $\tSpan\{\bar{X}_i,\bar{T}_i\bar{X}_i\}$ with the generalized eigenspaces are $d_1,\ldots,d_i-1,\ldots,d_{r+1}.$ If $v_i\notin X,$ then $\tSpan\{X,TX\}=\tSpan\{X,v_i\}$ is $T$-stable, violating Condition \eqref{eq:XisnotTstable}. Hence $v_i\in X$ and $\bar{X}_i$ is $n-1$ dimensional. We denote this reduction step by $$L'^{f,T,*}_{\{d_1,\ldots,d_{r+1}\}}(k^s) \xrightarrow[\delta]{\,\,\sim\,\,} L'^{f/(x-\al_i)^2,\bar{T}_i,*}_{\{d_1,\ldots,d_i-1,\ldots,d_{r+1}\}}(k^s).$$
$\delta$ is bijective, its inverse is given by taking the pre-image of the projection map $H_i\rightarrow \bar{U}_i.$

The stabilizers are affected in the same manner as in the odd case. We summerize this reduction step in the following lemma.

\begin{lemma}\label{lem:d_reductioneven}
Suppose $d_i\geq1,$ then there is a bijection
$$L'^{f,T,*}_{\{d_1,\ldots,d_{r+1}\}}(k^s) \xrightarrow[\delta]{\,\,\sim\,\,} L'^{f/(x-\al_i)^2,\bar{T}_i,*}_{\{d_1,\ldots,d_i-1,\ldots,d_{r+1}\}}(k^s).$$ The sizes of the stabilizers do not change, unless $m_i=2,d_i=1$ in which case it decreases by a factor of 2.
\end{lemma}

This reduction can be described projectively as intersecting the quadric defined by $Q$ with the tangent plane to $v$, then projecting away from $v$. Such an operation does not preserve the rulings. Two (projective) $n$-planes in $Q$ lying in the same ruling could be sent to different rulings via this procedure. For example take a smooth quadric in $\bbp^7,$ and two $3$-planes $Y_1,Y_2$ on the quadric intersecting at a line. Then these two $3$-planes lie on the same ruling. If the tangent plane to $v$ contains this line, then the images of $Y_1,Y_2$ lie in different rulings since their intersection codimension is 1. If the tangent plane to $v$ meets this line at a point, then the images $Y_1,Y_2$ lie in the same ruling as their intersection codimension is 2. Similar examples can be written down when $Y_1,Y_2$ lie on different rulings.

\subsubsection*{Reduction on $f$}

By the above reduction step, it remains to study $L'^{f,T,*}_{\{0,0,\ldots,0\}}(k^s).$ We will describe the reduction map, state the corresponding result, then give the proof. There is a slight difference to the odd case due to dimension reasons. Once again, the proof is just hardcore linear algebra, so we recommend the interested reader to prove it himself.

Suppose $\al$ is a root of $f$ of multiplicity $m\geq2.$ Let $X\in L'^{f,T,*}_{\{0,0,\ldots,0\}}(k^s)$ be arbitrary. Let $v$ denote an eigenvector of $T$ with eigenvalue $\al.$ Suppose $v'\in U$ such that $(T-\al)v'=v.$ Since $b(v,v)=0,$ we can consider the descent to $\bar{U}=v^\perp/v.$ As in the above reduction step, $Q$ descends to a non-degenerate quadratic form $\bar{Q}$ on $\bar{U}$ and $T$ descends to a regular self-adjoint operator $\bar{T}$ on $\bar{U}$ with characteristic polynomial $f(x)/(x-\al)^2.$

Observe that $v\notin \tSpan\{X,TX\}$ since $\tSpan\{X,TX\}$ contains no non-zero $T$-stable subspace. Therefore the map $U\rightarrow U/v$ is bijective when restricted to $\tSpan\{X,TX\}$. Denote the image of $X\cap v^\perp$ in $\bar{U}=v^\perp/v$ by $\bar{X}.$ As in the above reduction step, $\tSpan\{\bar{X},\bar{T}\bar{X}\}$ is an $n$-dimensional isotropic subspace of $\bar{U}$.

\begin{lemma}\label{lem:reductionwell}
$\tSpan\{\bar{X},\bar{T}\bar{X}\}$ has no non-zero $\bar{T}$-stable subspace.
\end{lemma}

\bpr Its only possible non-zero $\bar{T}$-stable subspace is the line spanned by $v'+<v>$. Suppose for a contradiction that $v'+cv\in\tSpan\{X,TX\}$ for some $c\in k$. Since $\tSpan\{X,TX\}$ has no non-zero $T$-stable subspace, we see that $v',v'+cv\notin X.$ Since $\tSpan\{X,TX\}$ is isotropic, we see that $v'+cv$ is orthogonal to every element in $(T-\al)X,$ and hence $v$ is orthogonal to every element of $X$. Since $v'+cv$ also lies in $X^\perp,$ we see that $v'\in X^\perp.$ Finally, $b(v,v')=0$ a priori due to the assumption that $v'+<v>\in \bar{U}.$ Combining these, one concludes that the $(n+2)$-dimensional subspace $\tSpan\{X,v',v\}$ is isotropic in $U$ with respect to $b$, contradicting to the fact that $U$ only has dimension $2n+2$.\epr

Consequently, $X\nsubseteq v^\perp$, for if otherwise $\bar{X}=\tSpan\{\bar{X},\bar{T}\bar{X}\}$ for dimension reasons and hence is $\bar{T}$-stable, which contradicts both Lemma \ref{lem:reductionwell} and Condition \ref{eq:XisnotTstable}. One now has the following well-defined map.

\begin{lemma}\label{lem:f_reductioneven}
Suppose $n\geq2$. The map sending $X$ to $\bar{X}$ defines a surjection
$$L'^{f,T,*}_{\{0,0,\ldots,0\}}(k^s) \xrightarrow{\quad} L'^{f/(x-\al)^2,\bar{T},*}_{\{0,0,\ldots,0\}}(k^s).$$
This map is bijective if $m>2$ and is two-to-one if $m=2.$ In both cases, $$|\tStab_{J_T}(X)| = |\tStab_{J_{\bar{T}}}(\bar{X})|,\qquad \mbox{for any } X\in L'^{f,T}_{\{0,\ldots,0\}}(k^s).$$
\end{lemma}

\bpr We first prove surjectivity. Suppose $\bar{X}\in L'^{f/(x-\al)^2,\bar{T},*}_{\{0,0,\ldots,0\}}(k^s).$ Let $b_\al$ denote the bilinear form $$b_\al(u,u')=b(u,(T - \al)u').$$
Since $v$ lies in the kernel of $b_\al$, we see that $b_\al$ descends to a non-degenerate bilinear form on the $2n+1$ dimensional vector space $U/v.$ Denote by $\perp_\al$ the perpendicular space with respect to $b_\al.$ Suppose for a contradiction that $\tSpan\{\bar{X},\bar{T}\bar{X}\}$ is isotropic with respect to $b_\al$. Then inside $\bar{U},$ $$\bar{T}^2\bar{X}\subset\tSpan\{\bar{X},\bar{T}\bar{X}\}^{\perp}=\tSpan\{\bar{X},\bar{T}\bar{X}\}.$$
Hence the entire $\tSpan\{\bar{X},\bar{T}\bar{X}\}$ is $\bar{T}$-stable. Contradiction.

Observe that $b_\al$ descends to a non-degenerate bilinear form on the $2$-dimensional vector space $\bar{Y}=\tSpan\{\bar{X},\bar{T}\bar{X}\}^{\perp_\al}/\bar{X}.$ Indeed a priori, $b_\al$ descends to a non-degenerate form on $\bar{X}^{\perp_\al}/\bar{X},$ and $\bar{X}^{\perp_\al}$ is spanned by $\tSpan\{\bar{X},\bar{T}\bar{X}\}^{\perp_\al}$ and a non-isotropic vector $u$ in $\bar{T}\bar{X}.$ Given any $w\in\tSpan\{\bar{X},\bar{T}\bar{X}\}^{\perp_\al},$ one can first find a $w'\in\bar{X}^{\perp_\al}$ such that $b_\al(w,w')\neq0,$ then adjust $w'$ by a multiple of $u$ so it lands in $\tSpan\{\bar{X},\bar{T}\bar{X}\}.$

As a 2-dimensional non-degenerate quadratic space, $\bar{Y}$ has two 1-dimensional isotropic lines, denote by $\bar{X}_1,\bar{X}_2$ their pre-images in $\tSpan\{\bar{X},\bar{T}\bar{X}\}^{\perp_\al}.$

Suppose $m\geq3,$ then as in the odd case, $b_\al(v',v')=b(v',v)=0,$ so up to renaming, $\bar{X}_1=\tSpan\{v'+<v>,\bar{X}\}\subset v^\perp/v.$ Since $\tSpan\{\bar{X}_1,\bar{X}_2\}$ has dimension $n+1,$ it is not isotropic with respect to $b_\al.$ Therefore, $b_\al(w,v')=b(w,v)\neq0$ for some $w+<v>\in\bar{X}_2.$ Up to scaling, we may assume $b(w,v)=1$ and by replacing $w$ by $w-\frac{1}{2}b(w,w)v,$ we may also assume $b(w,w)=0.$ Consider
\begin{eqnarray*}
X^w &=& \tSpan\{w,u - b(w,u)v\}_{u + <v>\in \bar{X}}\subset U,\\
(T-\al)X^w &=& \tSpan\{(T - \al)w, (T - \al)v\}.
\end{eqnarray*}
It is clear that $\tSpan\{X^w,TX^w\}$ is isotropic with respect to $b$ by the construction of $w$. Since $w\notin v^\perp,$ we have $\bar{X^w}=\bar{X}.$ Since $b(w, c_2v)=c_2$, we see that $\tSpan\{X^w,TX^w\}$ contains no elements of the form $c_2v$ since it is isotropic. Therefore $\tSpan\{X^w,TX^w\}$ has no non-zero $T$-stable subspace. We have now proved surjectivity when $m\geq3.$

Suppose now $X'\in L^{f,T,*}_{\{0,\ldots,0\}}(K)$ maps to $\bar{X}$. Then the image of $X'$ in $U/v$, denoted suggestively by $\bar{X'}_2$ is an $n$-plane isotropic to $b_\al,$ it contains $\bar{X}$ and is $b_\al$-orthogonal to $\tSpan\{\bar{X},\bar{T}\bar{X}\}.$ Since it does not contain $v'+<v>,$ we conclude that $\bar{X'}_2=\bar{X}_2.$ Since the process from $\bar{X}_2$ to $X^w$ is just adjusting with the correct multiples of $v$, we see that $X'=X^w.$ The way how the stabilizer changes is identical to the odd case.

We now deal with the case $m=2$. Write $\bar{X}_1=\tSpan\{w_1+<v>,\bar{X}\}$ and $\bar{X}_2 = \tSpan\{w_2+<v>, \bar{X}\}.$ We claim $w_1\notin v^\perp$ and likewise same with $w_2.$ Suppose for a contradiction that $w_1\in v^\perp.$ Since $\tSpan\{\bar{X},\bar{T}\bar{X}\}$ is not isotropic with respect to $b_\al,$ we see that $\tSpan\{\bar{X},\bar{T}\bar{X},w_1+<v>\}$ is an $n+1$ dimensional subspace of $v^\perp/v.$ As in the odd case, $b_\al$ is non-degenerate on $v^\perp/v$ because $T - \al$ acts invertibly on $v^\perp/v$. However, taking $\perp_\al$ inside $v^\perp/v,$ we see that $$\tSpan\{\bar{X},\bar{T}\bar{X},w_1+<v>\}^{\perp_\al}\supset\bar{X}_1.$$ The left hand side has dimension $n-1$ while the right hand side has dimension $n$. Contradiction.

Finally, we lift each $\bar{X}_i$ to $X^{w_i}$ by adding an appropriate multiples of $v.$ The resulting $X^{w_i}$ both maps to $\bar{X}$ under the reduction map. They are different from each other since their images in $U/v$ are different. Therefore we have proved surjectivity. The same argument as the above shows that $X^{w_1}$ and $X^{w_2}$ are precisely the two pre-images of $\bar{X}.$ Stabilizers behave in the same way as the odd case. \epr

\begin{cor}\label{cor:Lg0even}
$|L'^{f,T,*}_{\{0,0,\ldots,0\}}(k^s)| = 2^{r+1}$ and every element has trivial stabilizer in $J_T.$
\end{cor}

\bpr Apply the reduction steps like in the odd case. There are now five base cases which we illustrate as examples.\epr

\begin{example}\label{ex:1111}(Generic case)
Suppose reduction terminates with $f(x)=\prod_{i=1}^{r+1}(x-\al_i)$ with $r\geq3.$ In this case, one can apply the theory for the nonsingular case discussed in Example \ref{ex:nonWeier} and get $|L'^{f,T,*}|=2|L'^{f,T}|=2^{r+1}.$
\end{example}

\begin{example}\label{ex:211}
Suppose reduction terminates with $f(x) = (x - \al)(x - \be)(x - \gamma)^2.$ If one tries to apply reduction again on $\gamma,$ then $\bar{X}$ becomes 0-dimensional. Let $u,v,w_1$ denote the eigenvectors of $T$ with eigenvalue $\al,\be,\gamma$ respectively and let $w_2$ be such that $(T - \gamma)w_2=w_1.$ We seek coefficients $c_1,\ldots,c_4$ such that $X = <c_1u+c_2v+c_3w_1+c_4w_2>$ lies in $L'^{f,T,*}_{0,0,0}(k^s).$ Set
$$\Omega_1=b(u,u)\neq0,\quad\Omega_2=b(v,v)\neq0,\quad\Gamma_3=b(w_1,w_2)\neq0,\quad\Gamma_4=b(w_2,w_2).$$
Then the condition that $\tSpan\{X,TX\}$ is an isotropic $2$-plane becomes:
\[
\begin{pmatrix}
\Omega_1 & \Omega_2 & \Gamma_3 & \Gamma_4 \\
(\gamma-\al)\Omega_1 & (\gamma-\be)\Omega_2 & 0 & \Gamma_3\\
(\gamma-\al)^2\Omega_1 & (\gamma-\be)^2\Omega_2 & 0 & 0
\end{pmatrix}
\begin{pmatrix} c_1^2\\ c_2^2 \\ 2c_3c_4 \\c_4^2
\end{pmatrix} =
\begin{pmatrix} 0\\0\\ 0\\0
\end{pmatrix}
\]
Since $\Gamma_3,\Omega_1,\Omega_2$ are nonzero, the above matrix has a 1-dimensional kernel. Moreover, if any one of $c_1,c_2,c_4$ is zero, then they are all zero and $X$ is of the form $<c_3w_1>$ which does not lie in $L'^{f,T,*}_{0,0,0}(k^s).$ Now, given non-zero $c_1,c_2,c_4,$ one gets a unique solution for $c_3.$ Therefore, there are $8=2^3$ choices for $X$ depending on which square roots one chooses for $c_1,c_2,c_4.$
\end{example}

\begin{example}\label{ex:31}
Suppose reduction terminates with $f(x) = (x - \al)^3(x - \be).$ Let $u_1,v$ denote the eigenvectors of $T$ with eigenvalue $\al,\be$ respectively and let $u_2,u_3$ be such that $(T-\al)^2u_3=(T-\al)u_2=u_1.$ We seek coefficients $c_1,\ldots,c_4$ such that $X = <c_1u_1+c_2u_2+c_3u_3+c_4v>$ lies in $L'^{f,T,*}_{0,0}(k^s).$ Set
$$\Omega=b(v,v)\neq0,\quad\Gamma_4=b(u_1,u_3)=b(u_2,u_2)\neq0,\quad\Gamma_5=b(u_2,u_3),\quad\Gamma_6=b(u_3,u_3).$$
Then the condition that $\tSpan\{X,TX\}$ is an isotropic $2$-plane becomes:
\[
\begin{pmatrix}
\Gamma_4 & \Gamma_5 & \Gamma_6 & \Omega \\
0 & \Gamma_4 & \Gamma_5 & (\be-\al)\Omega\\
0 & 0 & \Gamma_4 & (\be-\al)^2\Omega
\end{pmatrix}
\begin{pmatrix} c_2^2+2c_1c_3\\ 2c_2c_3 \\ c_3^2 \\c_4^2
\end{pmatrix} =
\begin{pmatrix} 0\\0\\ 0\\0
\end{pmatrix}
\]
Since $\Gamma_4,\Omega$ are nonzero, the above matrix has a 1-dimensional kernel and if any one of $c_2,c_3,c_4$ is zero, then all of them are zero and $X$ is of the form $<c_1u_1>$ which does not lie in $L'^{f,T,*}_{0,0}(k^s).$ Now, given non-zero $c_3,c_4,$ one gets a unique solution for $c_1,c_2.$ Therefore, there are $4=2^2$ choices for $X$ depending on which square roots one chooses for $c_3,c_4.$
\end{example}

\begin{example}\label{ex:22}
Suppose reduction terminates with $f(x) = (x - \al)^2(x - \be)^2.$ Let $u_1,v_1$ denote the eigenvectors of $T$ with eigenvalue $\al,\be$ respectively and let $u_2,v_2$ be such that $(T-\al)u_2=u_1,(T-\be)v_2=v_1.$ We seek coefficients $c_1,\ldots,c_4$ such that $X = <c_1u_1+c_2u_2+c_3v_1+c_4v_2>$ lies in $L'^{f,T,*}_{0,0}(k^s).$ Set
$$\Gamma_3=b(u_1,u_2)\neq0,\quad\Gamma_4=b(u_2,u_2),\quad\Omega_3=b(v_1,v_2)\neq0,\quad\Omega_4=b(v_2,v_2).$$
Then the condition that $\tSpan\{X,TX\}$ is an isotropic $2$-plane becomes:
\[
\begin{pmatrix}
\Gamma_3 & \Gamma_4 & \Omega_3 & \Omega_4 \\
0 & \Gamma_3 & (\be-\al)\Omega_3 & \Omega_3+(\be-\al)\Omega_4\\
0 & 0 & (\be-\al)^2\Omega_3 & 2(\be-\al)\Omega_3+(\be-\al)^2\Omega_4
\end{pmatrix}
\begin{pmatrix} 2c_1c_2\\ c_2^2 \\ 2c_3c_4 \\c_4^2
\end{pmatrix} =
\begin{pmatrix} 0\\0\\ 0\\0
\end{pmatrix}
\]
Since $\Gamma_3,\Omega_3$ are nonzero, the above matrix has a 1-dimensional kernel. If any one of $c_2,c_4$ is zero, then both of them are zero and $X$ is of the form $<c_1u_1+c_3v_1>$. In this case, $\tSpan\{X,TX\}$ either contains $<u_1>$ or $<v_1>$ both of which are $T$-stable thereby forcing $X\notin L'^{f,T,*}_{0,0}(k^s).$ Note if $c_1$ and $c_3$ are both non-zero, then $X$ satisfy the weaker condition that $X$ contains no non-zero $T$-stable subspace. Moreover, $X\in L^{f,T,*}_{1,1}(k^s)$ violates Condition \eqref{eq:XisnotTstable}. It is clear that there are infinitely many such $X$.

Now, given non-zero $c_2,c_4,$ one gets a unique solution for $c_1,c_3.$ Therefore, there are $4=2^2$ choices for $X$ depending on which square roots one chooses for $c_2,c_4.$
\end{example}

\begin{example}\label{ex:400}
Suppose reduction terminates with $f(x) = (x - \al)^4$ Let $u_1$ denote the eigenvector of $T$ with eigenvalue $\al$ and let $u_2,u_3,u_4$ be such that $$(T-\al)^3u_4=(T-\al)^2u_3=(T-\al)u_2=u_1.$$ We seek coefficients $c_1,\ldots,c_4$ such that $X = <c_1u_1+c_2u_2+c_3u_3+c_4u_4>$ lies in $L^{f,T,*}_{0}(K).$ Set
$$\Gamma_5=b(u_1,u_4)\neq0,\quad\Gamma_6=b(u_2,u_4)=b(u_3,u_3),\quad\Gamma_7=b(u_3,u_4),\quad\Gamma_8=b(u_4,u_4).$$
Then the condition that $\tSpan\{X,TX\}$ is an isotropic $2$-plane becomes:
\[
\begin{pmatrix}
\Gamma_5 & \Gamma_6 & \Gamma_7 & \Gamma_8 \\
0 & \Gamma_5 & \Gamma_6 & \Gamma_7\\
0 & 0 & \Gamma_5 & \Gamma_6
\end{pmatrix}
\begin{pmatrix} 2c_1c_4+2c_2c_3\\ 2c_2c_4+c_3^2 \\ 2c_3c_4 \\c_4^2
\end{pmatrix} =
\begin{pmatrix} 0\\0\\ 0\\0
\end{pmatrix}
\]
Since $\Gamma_5$ is nonzero, the above matrix has a 1-dimensional kernel and if $c_4$ is zero, then $c_3$ is also zero. In this case, any $X$ of the form $<c_1u_1+c_2u_2>$ solves the above equation. However, for all such lines, $\tSpan\{X,TX\}$ contains the $T$-stable subspace $<u_1>$ thereby forcing $X\notin L'^{f,T,*}_{0}(k^s).$ Note if $c_1\neq0,$ then $\tSpan\{X,TX\}=<u_1,u_2>$ and $X\in L'^{f,T,*}_{2}(k^s)$ violating Condition \ref{eq:XisnotTstable} while all such $X$ still satisfy the weaker condition that it contains no non-zero $T$-stable subspace.

Now, given a non-zero $c_4,$ one gets a unique solution for $c_1,c_2,c_3.$ Therefore, there are $2=2^1$ choices for $X$ depending on which square root one chooses for $c_4.$
\end{example}

\textbf{Proof of Theorem \ref{thm:regularcooleven}: }Applying Lemma \ref{lem:d_reductioneven} repeatedly gives a bijection
$$L'^{f,T,*}_{\{d_1,\ldots,d_{r+1}\}}(k^s)\xrightarrow[\delta]{\,\,\sim\,\,} L'^{\prod_i(x-\al_i)^{m_i-2d_i},T',*}_{\{0,0,\ldots,0\}}(k^s),$$ and for any $X\in L'^{f,T,*}_{\{d_1,\ldots,d_{r+1}\}}(k^s),$ $$|\tStab_{J_T}(X)| = 2^a|\tStab_{J_{T'}}(\delta(X))|.$$

The polynomial $g(x)=\prod_i(x-\al_i)^{m_i-2d_i}$ has $r+1-a$ distinct roots, hence applying Corollary \ref{cor:Lg0even} to $g$ then dividing by 2 to go from $|L'^{f,T,*}|$ to $|L'^{f,T}|$ completes the proof.\epr
\end{example}